\documentclass[reqno]{amsart}
\usepackage{amsmath}
 \usepackage{amssymb}
\newcommand{\mysection}[1]{\section{#1}
      \setcounter{equation}{0}}

\newcommand{\nliminf}{\operatornamewithlimits{\underline{lim}}}

\newtheorem{theorem}{Theorem}[section]
\newtheorem{lemma}[theorem]{Lemma}

\newtheorem{corollary}[theorem]{Corollary} 

\theoremstyle{definition}
\newtheorem{assumption}{Assumption}[section]

\theoremstyle{remark}
\newtheorem{remark}{Remark}[section]

\newtheorem{example}{Example}[section]

\newcommand{\tr}{\text{\rm tr}\,}
\newcommand{\loc}{\text{\rm loc}}
\newcommand{\sign}{\text{\rm sign}}
\newcommand\esssup{\operatornamewithlimits{ess\,sup\,}}

 \makeatletter 
 \def\dashint{%  
 \operatorname%
 {\,\,\text{\bf--}\kern-.98em\DOTSI\intop\ilimits@\!\!}}
\def\ninf{\qopname\relax\@empty{inf\phantom{p}\!\!\!}}

 \makeatother

\def\dashnorm{\,\,\text{\bf--}\kern-.5em\|}

\newcommand\bR{\mathbb{R}}

\newcommand\bS{\mathbb{S}}

\newcommand\bZ{\mathbb{Z}}

\newcommand\cA{\mathcal{A}}

\def\cC{\mathcal{C}}

\newcommand\cF{\mathcal{F}}

\newcommand\cZ{\mathcal{Z}}

\newcommand\frB{\mathfrak{B}}

\begin{document}

\title[{Nondegenerate It\^o
 processes with moderated drift}]
{On  nondegenerate It\^o
 processes with moderated drift}

\author{N.V. Krylov}
 
\email{nkrylov@umn.edu}
\address{127 Vincent Hall, University of Minnesota,
 Minneapolis, MN, 55455}
 


\begin{abstract}
In this paper we present an approach to
proving parabolic Aleksandrov estimates
with mixed norms
for stochastic integrals with singular
``moderated'' drift.
\end{abstract}

\maketitle

\mysection{Introduction}

Let  $d_{1}$ be an integer $\geq2$,
$(\Omega,\cF,P)$ be a complete probability space,
and let $(w_{t},\cF_{t})$ be a $d_{1}$-dimensional
Wiener process on this space with complete, relative to
$\cF,P$, $\sigma$-fields $\cF_{t}$. Fix $\delta\in(0,1]$
and denote by $\bS_{\delta}$ the set of symmetric $d\times d$-matrices whose eigenvalues are in $[\delta,\delta^{-1}]$. Let $\sigma_{t},t\geq0$,
be a progressively measurable process with values in the set 
of $d\times d_{1}$-matrices such that
$a_{t}:=(1/2)\sigma_{t}\sigma^{*}_{t}\in \bS_{\delta}$ for all $(\omega,t)$, and let $b_{t},t\geq0$, be an $\bR^{d}$-valued
progressively  measurable process.
Assume that for any $T\in[0,\infty)$ and $\omega$
\begin{equation}
                                      \label{1.7.2}
\int_{0}^{T} |b_{t}| \,dt<\infty.
\end{equation}

Under these conditions the stochastic process  
$$
x_{t}=\int_{0}^{t}\sigma_{s}\,dw_{s}
+\int_{0}^{t}b_{s}\,ds
$$
is well defined. It is well known that, if $b$ is bounded,
then for any $T\in(0,\infty)$ and Borel $g(x)\geq0$
given on $\bR^{d}$ and
$f(t,x)\geq0$ given on $\bR^{d+1}=\{(t,x):t\in\bR,x\in\bR^{d}\}$
\begin{equation}
                                        \label{1.7.1}
E\int_{0}^{T}f(x_{t})\,dt\leq N\|f\|_{L_{d} },\quad E\int_{0}^{T}f(t,x_{t})\,dt\leq N\|f\|_{L_{d+1} },
\end{equation}
where the constant $N$ is independent of $g$ and $f$
and we use the same symbol $L_{p}$
for $L_{p}(\bR^{d})$ and $L_{p}(\bR^{d+1})$.
Which $L_{p}$ is meant in each particular case
is easily determined by the number of
arguments of considered functions.

   There was a
considerable interest in reducing the $L_{d}$-norm of $g$
and the $L_{d+1}$-norm of $f$
  to the $L_{d_{0}}$-norm and the $L_{d_{0}+1}$-norm, respectively,
with
$d_{0}<d$, especially in the case of
stochastic equations and not just stochastic integrals. 
The main contributors to this line of research
in time homogeneous case (the coefficients depend only on $x_{s}$) are Cabr\'e
\cite{Ca_95} for bounded $b$,
  Fok   \cite{Fo_98} for
  $b\in L_{d+\varepsilon}$, the author
 in \cite{Kr_19_1}  
allowed $b\in L_{d}$, Hongjie Dong and the author
\cite{DK_21} allowed $b$ from a Morrey class.

The estimates like \eqref{1.7.1} are
indispensable in the theory of
controlled diffusion processes (see, for instance,
\cite{Kr_77}). First such  result in the time inhomogeneous
case with bounded $b$ and
$d$ (not $d_{0}$)  was published in \cite{Kr_74}. 
It  was extended by A.I. Nazarov and N.N. Ural'tseva
\cite{NU_85} to allow $b\in L_{d+1}$. The author in
\cite{Kr_86} developed a general approach to such estimates
and slightly improved the result of \cite{NU_85}.
By using this approach A.I. Nazarov in
\cite{Na_15}
for the first time proved estimates
in $L_{p,q}$-spaces with $b\in L_{p ,q }$.

 In \cite{CKS_00} by extending some earlier
results by Wang the authors, basically, prove the second
estimate in \eqref{1.7.1} in case $b$ is bounded with $d_{0}$ ($<d$) in place of
$d$. In \cite{4} it is extended to mixed norms of $f$
and $b\in L_{p,q}$ with $p,q$ smaller than in \cite{Na_15}.

 Here we present what should have been the start
of obtaining the results like in \cite{4} and in other papers
revolving around estimates like \eqref{1.7.1}.

For random variables $\xi,\eta$ we write $\xi\leq \eta$
if this inequality holds (a.s.) and use the notation
$$
E_{\cF_{t}}\eta=E\{\eta\mid \cF_{t}\},\quad
P_{\cF_{t}}(A)=E_{\cF_{t}}I_{A},\quad \partial_{t}=\frac{\partial}{\partial t},
$$
$$
D_{i}=\frac{\partial}{\partial x^{i}},\quad
D_{ij}=D_{i}D_{j},\quad Du=(D_{i}u),\quad D^{2}u=
(D_{ij}u).
$$

For simplicity we write simply $E_{\cF_{t}}\eta \leq A$ dropping $(a.s.)$.
Also we use the symbol $N$, sometimes with indices,
to denote various different constants and in the proofs
of our results they are supposed to be depending only on what is asserted in the particular result.

\mysection{Main assumption}
                                     \label{section 1.7.2}

For $t,R \in(0,\infty)$ and $x\in\bR^{d}$
denote by $\theta_{t}\tau'_{R}(x )$ the first exit time of
the process
$ x_{t+s}- x_{t}-x$, $s\geq0$, from $B_{R}=\{x\in\bR^{d}:|x|<R\}$ or,
equivalently, the first exit time of $x_{t+s}-x_{t} $
from $B_{R}( x)=B_{R}+x$.
Next, we set 
$\theta_{t}\tau_{R}(x)=R^{2}\wedge \theta_{t}\tau'_{R}(x)$.

If we are given a stopping time $\tau$, we set $(\theta_{\tau}\tau'_{R}(x)(\omega)=I_{\tau<\infty}\theta_{\tau(\omega)}\tau'_{R}(x)$. Similar sense is given to $\theta_{\tau}\tau_{R}(x)(\omega)$. 
If $x=0$, we drop it in the above notation.
It is useful to note that $t+\theta_{t}\tau'_{R}(x)$
are stopping times and
\begin{equation}
                                      \label{10.17.5}
\theta_{t}\tau'_{R }(x )\leq \theta_{t}\tau'_{2R }.
\end{equation}
Indeed, if $|x |\geq R$, the left-hand side is zero,
and if $|x|<R$, $B_R( x)\subset B_{2R}$.

Define one of the main quantities we will be using
by
$$
\bar b_{R}=\sup_{\rho\leq R}\hat b_{\rho},\quad
\hat b_{\rho}:=\frac{1}{\rho}\sup_{x\in\bR^{d}} 
\sup_{ t\geq0}
\esssup E_{\cF_{t}}\int_{0}^{\theta_{t}\tau_{ \rho}(x)}
|b_{t+s} |\,ds  .
$$
One can say that $\bar b_{R}$ ``moderates'' $b$ on scale $R$.
In light of \eqref{10.17.5}
\begin{equation}
                                      \label{10.17.6}
\bar b_{R}\leq \sup_{ \rho\leq R }
\frac{1}{\rho}\sup_{ t\geq0}
\esssup E_{\cF_{t}}\int_{0}^{\theta_{t}\tau_{ 2\rho }}
|b_{t+s} |\,ds  .
\end{equation}

\begin{remark}
                              \label{remark 10.14.1}
For any $\rho>0$ and stopping time $\tau$ 
and $\cF_{\tau}$-measurable $\bR^{d}$-valued $y$ we have
$$
E_{\cF_{\tau}}\int_{0}^{\theta_{\tau}\tau_{ \rho}(y)}
|b_{\tau+s} |\,ds
\leq \bar b_{\rho}\rho
$$
or, in other words, for any $A\in \cF_{\tau}$
\begin{equation}
                                   \label{10.14.1}
EI_{A}\int_{\tau}^{\tau+\theta_{\tau}\tau_{ \rho}(y)}
|b_{ s} |\,ds
\leq \bar b_{\rho}\rho P(A).
\end{equation}

Indeed, if $\tau$ and $y$ take only countably many values
(including $\infty$ for $\tau$),
estimate \eqref{10.14.1} immediately follows from the definition of $\bar b_{\rho}$. 
In the case of general $\tau$, one knows that there exists
a sequence of stopping times $\tau_{n}$
with values in the set of dyadic rationals such that
$\tau_{n}\downarrow \tau$. Furthermore, as is easy
to see, for any $\varepsilon\in(0,\rho)$
$$
\nliminf_{n\to\infty}\theta_{\tau_{n}}\tau_{ \rho}(y) 
\geq \theta_{\tau }\tau_{ \rho-\varepsilon}(y) 
, 
$$
$$ \int_{\tau}^{\tau+\theta_{\tau}\tau_{ \rho-\varepsilon}(y)} 
|b_{ s} |\,ds\leq \nliminf_{n\to\infty}
\int_{\tau_{n}}^{\tau_{n}+\theta_{\tau_{n}}\tau_{ \rho }(y)} 
|b_{ s} |\,ds
$$
and Fatou's lemma shows that \eqref{10.14.1} holds
if we replace $\rho$ in its left hand side with $\rho-\varepsilon$. After this replacement, still in the case of discrete
$y$, it will only remain
to use the monotone convergence theorem sending $\varepsilon\downarrow0$. To pass to general $y$ we approximate
it with the discrete ones $y_{n}$ and use the fact that
$$
\nliminf_{n\to\infty}\theta_{\tau}\tau_{\rho}(y_{n})
\geq \theta_{\tau }\tau_{ \rho-\varepsilon}(y) .
$$
 
\end{remark}

Next assumption, in which $m_{b}=m_{b}(d,\delta)\in(0,1]$ is a   number
to be specified later in Theorem \ref{theorem 8.2.1}, is  supposed to hold throughout
the paper after Theorem \ref{theorem 8.2.1}.
\begin{assumption}[$\bar b_{\rho_{b}}$]
                        \label{assumption 8.19.2}
  We have a $\rho_{b}\in(0,\infty)$ such that $\bar b_{\rho_{b}}
\leq m_{b}$.
\end{assumption}

\begin{remark}
                         \label{remark 8.21.01}
A very important feature of this assumption, which we however
are not going to use in this article, is
that it is preserved under self-similar dilations. 
\end{remark}

\begin{remark}
                                       \label{remark 3.27.1}
Usual way to deal with additive functionals shows 
 that
for any $n=1,2,...$, $\rho,t\in(0,\infty)$
$$
E_{\cF_{t}} \Big(\int_{0}^{\theta_{t}\tau_{ \rho }(x)}
|b_{t+s}|\,ds\Big)^{n} \leq n!\,\bar b_{\rho}^{n}\rho^{n}.
$$
Furthermore, by taking into account that
for any random variable $\xi\geq0$ and $\alpha
\in[1,2]$ we have $\big(E\xi^{\alpha}\big)^{1/\alpha}
\leq \big(E\xi \big)^{(2-\alpha)/\alpha}
\big(E\xi^{2}\big)^{(\alpha-1)/\alpha}$, we find that
for any $\varepsilon>0$ there exists $\alpha=\alpha(\varepsilon)
>1$ such that  
$$
\Big(E_{\cF_{t}} \Big(\int_{0}^{\theta_{t}\tau_{ \rho}(x)} 
|b_{t+s}|\,ds\Big)^{\alpha} 
 \Big)^{1/\alpha}\leq 2^{(\alpha-1)/\alpha}
\bar b_{\rho}\rho \leq (1+\varepsilon)\bar b_{\rho}\rho .
$$
 
\end{remark}

For $p,q\in[1,\infty)$  
 we introduce the space $L_{p,q} $ as the space of Borel
functions on $\bR^{d+1}=\bR\times \bR^{d}$ such that
$$
\|f\|^{q}_{L_{p,q} }:=
\int_{\bR}\Big(\int_{\bR^{d}} |f(t,x)|^{p}\,dx\Big)^{q/p}\,dt<\infty 
$$
if $p\geq q$ or
$$
\|f\|^{p}_{L_{p,q} }:=\int_ {\bR^{d}} \Big(\int_{\bR} 
|f(t,x)|^{q}\,dt\Big)^{p/q}\,dx<\infty 
$$
if $ p\leq q$ with natural interpretation
of these definitions if $p=\infty$ or $q=\infty$.
We write $f\in L_{p,q}(Q)$ to mean that $fI_{Q}\in L_{p,q}$.
Observe that
$p $ is associated with   $x$ and
$q$ with   $t$ and the interior
integral is always elevated to the power $\leq 1$. 
Naturally, $L_{p}(Q)=L_{p,p}(Q)$.

 Introduce
$$
  C_{ R}=[0,R^{2})\times B_{R},\quad
C_{ R}(t,x)=(t,x)+C_{ R}
$$
and let
$\cC_{R}$ be the collection of cylinders $C_{R}(t,x)$,
$(t,x)\in\bR^{d+1}$,  $\cC=
 \{\cC_{R}:R>0\} $.

Define
$$
\dashnorm h\|_{L_{p,q}(C)}=\|h\|_{L_{p,q}(C)}\|
1\|_{L_{p,q}(C)}^{-1}.
$$

\begin{lemma}
                       \label{lemma 12.11.1}
  Assume that 
there is a Borel function $h(t,x)$ 
such that $|b_{t}|\leq h(t,x_{t})$ and
for  some $\rho_{b}\in(0,\infty)$, any $\rho\leq \rho_{b}$,
and $C\in\cC_{\rho}$,  we have
\begin{equation}
                            \label{8.19.30}
\dashnorm h\|_{L_{d+1}(C)}\leq \bar h\rho^{-1},
\end{equation}
where $\bar h$ is a constant. Finally, let 
   $N(d,\delta )\bar h
\leq m_{b}$, where $N(d,\delta )$ is specified in the proof. Then Assumption   $(\bar  b_{\rho_{b}})$  \ref{assumption 8.19.2}
is satisfied.

\end{lemma}

Proof. By Theorem \ref{theorem 9.27.1} (below) for $\rho\leq\rho_{b}$
$$
I:=E_{\cF_{t}}\int_{0}^{\theta_{t}\tau_{ \rho}(x)}
|b_{t+s} |\,ds\leq 
E_{\cF_{t}}\int_{t}^{t+\theta_{t}\tau_{ \rho}(x)}
|h(s,x_{s}| \,ds
$$
$$
\leq NC\|h\|_{L_{d+1}(C_{\rho}(t,x_{t}+x))}
\leq NC\bar h\rho^{1/(d+1)},
$$
where
owing to the fact that the time-length of 
$C_{\rho}(t,x_{t}+x )$ is $\rho^{2}$, 
$$
C=\big(\rho^{2}+\hat b_{\rho}^{2}\rho^{2}\big)^{d/(2d+2)}\leq N
\rho^{d/(d+1)}(1+\hat b_{\rho})^{d/(d+1)},   
$$
which yields  $\hat b_{\rho}\leq N\bar h(1+\hat b_{\rho})^{d/(d+1)}$. This proves the lemma. \qed

In the last step of the above proof we tacitly
assumed that $\hat b_{\rho}<\infty$, which can always
be achieved by using stopping times associated with \eqref{1.7.2}.

\begin{remark}
                              \label{remark 12.11.1}

It turns out that \eqref{8.19.30} is satisfied with as small $\bar h$
as we like on the account of taking $|C|$
small if $\|h\|_{L_{d+2} }<\infty$. Indeed, by H\"older's
inequality, if $C\in\cC_{\rho}$, then
$$
\dashnorm h\|_{L_{d+1}(C)}\leq 
\dashnorm h\|_{L_{d+2}(C)}=
N(d)\rho^{-1}\|h\|_{L_{d+2}
(C) },
$$
where the last factor tends to zero as $\rho\to0$.
\end{remark}

On the other hand, it may happen that
\eqref{8.19.30} is satisfied  
but $h\not\in L_{d+2,\loc}$.

\begin{example}
                                          \label{example 8.24.1}
Take  $\alpha\in(0,d),\beta\in(0,1)$ such that
$\alpha+2\beta=d+1$ and consider the function
$g(t,x)=|t|^{-\beta}|x|^{-\alpha}$. Observe that
$$
\int_{C_{R}(t,x)}g(s,y)\,dyds  
=R  \int_{C_{1}(t',x')}g(s,y)\,dyds,
$$
where $t'=t/R^{2}$, $x'=x/R$. Obviously, the last integral
is a bounded function of $(t',x')$. Hence, the function
$h=g^{1/(d+1)}$ satisfies \eqref{8.19.30}. As is easy to see
for any $p>d+1$ one can find $\alpha$ and $\beta$
above such that $h\not\in L_{p,\loc}(\bR^{d+1})$.
 
\end{example}

However, there are very many situations when
\eqref{8.19.30} does not hold and 
Assumption   $(\bar  b_{\rho_{b}})$  \ref{assumption 8.19.2} is still satisfied.

\begin{example}
                    \label{example 5.23.1}

In $\bR^{d}$ with $d\geq 2$ take a $d$-dimensional
Wiener process $w_{t}$ and consider the system
$dx^{1}_{t}=dw^{1}_{t}+b(x_{t}) \,dt$,
$dx^{i}_{t}=dw^{i}_{t}$, $i\geq2$, where
$$
b(x )=\beta(x^{1} ),\quad \beta(r)=-|r|^{-\alpha}I_{(-1,1)}(r) \sign\, r,
$$ 
  and $\alpha<1$ is as close
to $1$ as we wish. The solutions of our system
form a strong Markov time-homogeneous process for which
Assumption   $(\bar  b_{\rho_{b}})$  \ref{assumption 8.19.2} is rewritten
as 
\begin{equation}
                                 \label{10.10.1}
\bar b_{\rho_{b}}:=\sup_{\substack{\rho\leq \rho_{b} 
\\ \,C\in\cC_{\rho}}}
\frac{1}{\rho}
\sup_{ x \in\bR^{d }}E_{ x}\int_{0}^{\tau_{C}}
|b( x_{s})|\,ds\leq m_{b},
\end{equation}
where $\tau_{C}$ is the first exit time of $(t,x_{t})$
from $C$.

Condition \eqref{8.19.30}
is definitely not satisfied if $\alpha$ is too close to $1$. 
However, if
$C=C_{\rho}(s,y)$ and $|y^{1} |\leq 2\rho$,
then with $P_{ x}$-probability one
$\tau_{C}\leq \tau_{(-3\rho,3\rho)}$, where by
$\tau_{(a,b)}$ we denote the first exit time of $x^{1}_{t}$ from $(a,b)$. In that case by using It\^o's formula one gets that, for $(0,x)\in C$
$$
 E_{ x}\int_{0}^{\tau_{C}}|b(x_{s})|\,ds\leq E_{ x}\int_{0}^{\tau_{(-3\rho,3\rho)}}
|b(x_{s})|\,ds=:\phi(|x^{1}|),
$$
where $\phi(r)=0$ if $|x^{1}|\geq 3\rho$ and  otherwise
$$
\phi(r)=r- 3\rho+ \int_{r}^{3\rho}
\exp\Big(\frac{2}{1-\alpha}(t\wedge 1)^{1-\alpha}\Big)
\,dt\leq N\rho.
$$

In case $ y^{1} \geq 2\rho$ observe that,
for $r_{\pm}:=y^{1}\pm \rho$,
$$
E_{ x}\int_{0}^{\tau_{C}}|b(x_{s})|\,ds\leq E_{ x}\int_{0}^{\tau_{( r_{-},r_{+})}}
|b(x_{s})|\,ds=: \psi(x^{1}) ,
$$
which is zero if $x^{1}\not \in ( r_{-},r_{+}) $
and if $x^{1}  \in ( r_{-},r_{+}) $
 by It\^o's formula equals (observe that $|\beta(r)|=
-\beta(r)$ on $( r_{-},r_{+}) $)
$$
\psi(x^{1})=-E_{t,x}(x^{1}_{\tau_{( r_{-},r_{+})}}-x^{1}_{0})
$$
which  is less than $2\rho$.
Therefore, $\psi(x^{1})\leq 2\rho$ and,
since similar situation occurs if $y^{1}\leq -2\rho$,
$\bar b_{\infty}<\infty$.

To show that \eqref{10.10.1} is satisfied 
we show that $\bar b_{\rho}\to0$ as $\rho\to0$.

Note that $\phi(r)$ is a decreasing function for
$r\geq0$ and if $3\rho\leq1$
$$
\phi(r)\leq\phi(0)=\int_{0}^{3\rho}\Big[
\exp\Big(\frac{2}{1-\alpha}(t\wedge 1)^{1-\alpha}\Big)
-1\Big]\,dt
$$
$$
\leq 3\rho\Big[
\exp\Big(\frac{2}{1-\alpha}(3\rho)^{1-\alpha}\Big)
-1\Big]\leq N\rho^{2-\alpha}.
$$

To estimate $\psi$ for $r \in ( r_{-},r_{+}) $
introduce
$$
 \xi(r)=r-r_{+}+\eta(r),\quad \eta(r)=\frac{2\rho}{  \gamma -1}
\big(e^{2\hat b(r-r_{+})}-1 \big),
$$
where
$$
\hat b= ( 2/|y^{1}| )^{\alpha},\quad \gamma=e^{-4\hat b\rho}.
$$
By observing that on $( r_{-},r_{+}) $ we have $\beta(r)\geq -\bar b$ and $\eta'\leq0$, we obtain
$$
(1/2)\eta''+b\eta'\leq (1/2)\eta''-\hat b \eta'=0,
\quad (1/2)\xi''+b\xi'\leq b=-|b|.
$$
Furthermore, $\xi(r_{\pm})=0$ and a simple application 
of It\^o's formula shows that
$\psi\leq \xi$ on $ ( r_{-},r_{+}) $.

To estimate $\xi$ use that $e^{t}-1\geq t$
implying that $\eta(r)\leq 4\rho\hat b(r_{+}-r)(1-\gamma)^{-1}$, so that
$$
\xi(r)\leq (r_{+}-r)\big(4\rho\hat b(1-\gamma)^{-1}-1\big)\leq 2\rho \big(4\rho\hat b(1-\gamma)^{-1}-1\big).
$$
Note that $\hat b\rho\leq \rho^{1-\alpha}$ since $y^{1}
\geq 2\rho$ and $\psi(r)\leq\xi(r)\leq 10\rho^{2-\alpha}$ if $\rho$
is small enough.
\end{example}

\mysection{Preliminary results}

                              \label{section 10.25.2}

We use the following which 
combines particular cases of conditional versions of  Lemmas 4.1 and 4.2 of \cite{Kr_20_2}.
We have two stopping times $\gamma\geq\tau$.

\begin{theorem}
                                        \label{theorem 9.27.1}

For any $\lambda\geq 0$ and  Borel $f(s,y),g(y)\geq0$  
\begin{equation}
                                   \label{5.6.401}
  E_{\cF_{\tau}}\int_{\tau}^{\gamma} e^{-\lambda(s-\tau)}  
f(  s,x_{ s})\,ds\leq N(d, \delta  )
\big(A_{\lambda} + B ^{2}_{\lambda}
\big)^{d/(2d+2 )}
\|f\|_{L_{d+1 }},
\end{equation}
\begin{equation}                      \label{3.7.1} 
  E_{\cF_{\tau}}\int_{\tau}^{\gamma}e^{-\lambda(s-\tau)}   
g(  x_{ s})\,ds\leq N(d, \delta  )
\big(A_{\lambda} + B_{\lambda}^{2}
\big)^{1/2}
\|g\|_{L_{d  } },
\end{equation}
where  
$$
A_{\lambda}=E_{\cF_{t}}\int_{\tau}^{\gamma}e^{-\lambda(s-\tau)}\,ds,\quad
B_{\lambda}=E_{\cF_{\tau}}\int_{\tau}^{\gamma}e^{-\lambda(s-\tau)}|b_{ s}|\,ds.
$$

\end{theorem}

\begin{remark}
                           \label{remark 9.23.1}
It turns out that, if $\gamma-\tau\leq \theta_{\tau}\tau'_{ R}$
($\theta_{\tau}\tau'_{ R}$ is the first exit time of $x_{ \tau +s}$
from $B_{R}(x_{\tau})$), then $2d\delta A_{0}\leq R^{2}+2RB_{0}\leq  B_{0}^{2}+2R^{2}$,
so that for $\lambda=0$, $A_{0}$ in \eqref{5.6.401} and \eqref{3.7.1}   can be replaced
with $R^{2}$. 

Indeed, by It\^o's formula
$$
|x_{t\wedge\gamma}-x_{t\wedge\tau}|^{2} = 2\int_{t\wedge\tau}^{t\wedge\gamma}
 \tr a_{s} \,ds
+2
\int_{t\wedge\tau}^{t\wedge\gamma}
 (x^{i}_{s}-x^{i}_{t\wedge\tau})b^{i}_{s}\,ds
+2\int_{t\wedge\tau}^{t\wedge\gamma} 
(x^{i}_{s}-x^{i}_{t\wedge\tau})\sigma^{ik}_{s}\,dw^{k}_{s}.
$$
Here the stochastic integral is a local martingale.
Therefore, by replacing $t$ with $ \tau_{n}$ for an appropriate sequence
of stopping times  $\tau_{n}\to\infty$, then taking the expectations and using that
$$
|x_{\tau_{n}\wedge\gamma}-x_{\tau_{n}\wedge\tau}|^{2} \leq  R^{2},
$$
$$
E_{\cF_{\tau}}\Big|\int_{\tau_{n}\wedge\tau}^{\tau_{n}\wedge\gamma}
 (x^{i}_{s}-x^{i}_{\tau_{n}\wedge\tau})b^{i}_{s}\,ds\Big|
=I_{\tau_{n}\geq\tau}E_{\cF_{\tau}}\Big|\int_{\tau_{n}\wedge\tau}^{\tau_{n}\wedge\gamma}
 (x^{i}_{s}-x^{i}_{\tau_{n}\wedge\tau})b^{i}_{s}\,ds\Big|
\leq RB_{0},
$$
we find
$$
2E_{\cF_{\tau}}\int_{\tau_{n}\wedge\tau}^{\tau_{n}\wedge\gamma}
  \tr a_{s} \,ds
\leq R^{2}+2RB_{0}.
$$
Sending $n\to\infty$ yields our claim.
\end{remark}
 
Observe that if $\gamma=\tau+\theta_{\tau}\tau_{ R}(y) $ 
in Theorem \ref{theorem 9.27.1}, 
then obviously $\gamma-\tau
\leq R^{2}$ and $A_{0} \leq R^{2}$. In that case also
$B_{0} \leq \bar b_{R}R$ by definition. Hence
we have the following.

\begin{lemma}
                                      \label{lemma 8.16.1}
 
For any Borel $f,g \geq0$, 
$\cF_{\tau}$-measurable $\bR^{d}$-valued $y$, and $R>0$ 
we have
\begin{equation}
                                          \label{9.29.2}
 E_{\cF_{\tau}}\int_{0}^{\theta_{\tau}\tau_{ R}(y)  }  
f( \tau+ s,x_{\tau+ s})\,ds\leq 
 N(d,  \delta)(1+\bar b_{R}) ^{d/  (d+1)  }
 R 
 ^{d/ (d+1)  }
\|f\|_{L_{d+1 }},
\end{equation}
\begin{equation}
                                            \label{8.22.2}
  E_{\cF_{\tau}}\int_{0}^{\theta_{\tau}\tau_{ R}(y)  }   
g(  x_{\tau+s})\,ds\leq N(d, \delta  )
(1+\bar b_{R})  
 R 
\|g\|_{L_{d  } },
\end{equation}
 
\end{lemma}
 
 Obviously, \eqref{9.29.2} and \eqref{8.22.2}
hold also if $f$ and $g$ are random,
provided that they are jointly measurable and $\cF_{\tau}$-measurable for other arguments fixed. By taking
$f(s,x)=I_{C_{R}(\tau,x_{\tau}+y)}(s,x)$ we get from \eqref{9.29.2} that
\begin{equation}
                                          \label{3.29.1}
E_{\cF_{\tau}}\theta_{\tau}\tau_{ R }(y)\leq N(d,\delta)(1+\bar b_{R})^{d/(d+1)}R^{2}.
\end{equation}

Estimate \eqref{3.29.1} says that   $\theta_{\tau}\tau_{ R }$ is of order
not more than $R^{2}$ for small $R$.
A very important fact which is implied by Corollary
\ref{corollary 7.29.1} is that 
$\theta_{\tau}\tau_{R}$ is of order
not less than $R^{2}$.
To show this we need  the following result,
in which
\begin{equation}
                                                  \label{1.6.1}
\theta_{t}\gamma'_{ R}(x)  =\inf\{s\geq0:x_{t+s}-x_{t}   \in \bar B_{R}(x)  \}.
\end{equation}

\begin{theorem}
                           \label{theorem 8.2.1}
There are
   constants $\bar \xi=\bar \xi(d,\delta)\in (0,1) $ and
 $m_{b}=m_{b}(d,  \delta)\leq 1$ {\em
continuously\/} depending on $\delta$
such that if, for an $  R\in(0,\infty)$, we have
\begin{equation}
                                     \label{12.18.3}
  \bar b_{  R}\leq m_{b},
\end{equation}
then  for any stopping time $\tau$
and $\cF_{\tau}$-measurable $\bR^{d}$-valued $y$
\begin{equation}
                                          \label{8.2.2} 
  P_{\cF_{\tau}}( \theta_{\tau}\tau'_{R}(y)   \geq   R^{2} )\leq 1-\bar\xi,\quad
   P_{\cF_{\tau}}( \theta_{\tau}\tau'_{R}  \geq   R^{2} )\geq 
 \bar\xi I_{\tau<\infty}.   
\end{equation}
Moreover for $n=1,2,...$  
\begin{equation}
                                          \label{1.3.1} 
P_{\cF_{\tau}}( \theta_{\tau}\tau'_{ R}(y)  > nR^{2})
 \leq (1-\bar\xi)^{n},   
\end{equation}
so that 
\begin{equation}
                                          \label{3.7.2}
E _{\cF_{\tau}}\theta_{\tau}\tau'_{R}(y)  \leq N(d,\delta)R^{2},
\end{equation}
and
\begin{equation}
                                          \label{1.3.3}  
I:=E_{\cF_{\tau}}\int_{0}^{\theta_{\tau}\tau'_{ R}(y)}|b_{\tau+s}|\,ds
\leq N(d ,\delta) \bar b_{R}R.
\end{equation}

Furthermore,   on the set
 $\{|y|\leq 9R/16\}$  
\begin{equation}
                                          \label{1.2.1} 
P_{\cF_{\tau}}\big(\theta_{\tau}\tau'_{R}(y) \geq\theta_{\tau}\gamma'_{R/16}(y)    \big)\geq\bar\xi.   
\end{equation}

\end{theorem}

 We need an auxiliary result, in which
$$
m_{t,s}=-\int_{t}^{t+s}\sigma_{r}\,dw_{r}.
$$

\begin{lemma}
                                       \label{lemma 1.2.1}
(i) There exists $\kappa=\kappa(d)>0$ such that
for
$$
\psi_{t}(s,y)=R^{-4}\big(R^{2}-4|y|^{2}\big)^{2}\phi_{t,s},\quad
\phi_{t,s}=\exp\int_{t}^{t+s}
\kappa R^{-2}\,\tr a_{r}\,dr
$$
the process $\{\psi_{t}(s, m_{t,s}),\cF_{t+s}\}$ is a local
 submartingale for any $t\geq0$.

(ii) Take a $\zeta\in C^{\infty}_{0}(\bR)$
such that it is even, nonnegative, and decreasing
on $(0,\infty)$.
For  $T\in(0,\infty)$ and $x\in \bR$ and $t\leq T$define
$u(t,x)=E\zeta(x+w^{1}_{T-t}) )$. Also 
take $t\geq 0$, $x\in\bR^{d}$ and set
$$
\xi_{t,s}=\frac{(x+ m_{t,s},a_{t+s}(x+ m_{t,s}) )}
{|x+m_{s}|^{2}}\quad (0/0:=1),\quad
\eta_{t,s}=2\int_{0}^{s}\xi_{t,r}\,dr.
$$
Then   the process $\{u(\eta_{t,s},
|x+m_{t,s}|),\cF_{t+s}\}$
is a  supermartingale before $\eta_{t,s}$ reaches $T$,
in particular, on $[0,\delta^{2} T]$.

(iii) There exists $\alpha=\alpha(d,\gamma)>1$
such that for $u(x)=|x|^{-\alpha}$
and any $a\in\bS_{\gamma}$ we have
$
a^{ij}D_{ij}u(x)\geq0 
$ if $x\ne0$.

\end{lemma}

Proof. (i) It is easy to see that for a $\kappa=\kappa(d)>0$
we have $\kappa\mu^{2}-16\mu+32d^{-1}(1-\mu)\geq0
$ for all $\mu$, which implies that
for all $\lambda$
$$
\kappa(1-4\lambda^{2})^{2}-16(1-4\lambda^{2})
+128 d^{-1} \lambda^{2}\geq0.
$$

It follows that  (dropping $t$)
$$
R^{4}\phi_{s}^{-1}d\psi(s, m_{s} )=\kappa 
\big(R^{2}-4| m_{s}|^{2}\big)^{2}R^{-2}\tr a_{s}\,dt
$$
$$
-8\big(R^{2}-4|  m_{s}|^{2}\big)\big(2  m_{s}\,dm_{s}
+2\tr a_{s}\,ds\big)
+128( m_{s},a_{s}  m_{s})\,dt\geq dM_{s},
$$
where $M_{s}$ is a local martingale. This proves (i).

(ii) Observe that $u$ is smooth, even in $x$, and satisfies
$\partial_{t}u+(1/2)u''=0$.
Furthermore, as is easy to see $u'(t,x)\leq 0$
for $x\geq0$.
It follows by It\^o's formula that 
before $\eta_{t,s}$ reaches $T$ we have
(dropping obvious values of some arguments)
$$
du(\eta_{t,s},|x+m_{t,s}|)= \xi_{s} (2\partial_{t}u+u'')
\,ds +\frac{u'}{|x+m_{t,s}|} (\tr a_{t+s}-
 \xi_{s} )\,ds
+dM_{s},  
$$
where $M_{s}$ is a stochastic integral. Here the second
term with $ds$ is negative since $u'\leq0$, and this proves 
that $u(\eta_{t,s},|x+m_{t,s}|)$ is a local supermartingale.
Since it is nonnegative, it is a supermartingale.

Assertion (iii) is proved by simple computations.
The lemma is proved.\qed

{\bf Proof of Theorem \ref{theorem 8.2.1}}. 
While proving \eqref{8.2.2} we may concentrate
on the part of $\Omega$ where  $\tau $ is finite. Then  observe that for 
$$
\gamma:=R^{2}\wedge \inf\{s\geq0:|m_{\tau,s} |\geq R/2\}
$$
we have
$
\phi_{\gamma}\leq e^{\kappa d/\delta}.    
$
Hence, by Lemma \ref{lemma 1.2.1} (i) 
$$
1=\psi(0,0) \leq  E_{\cF_{\tau}}\psi(\gamma,m_{\tau,\gamma})
\leq e^{\kappa d/\delta^{2}}
P_{\cF_{\tau}}(\sup_{s\leq R^{2}}|m_{\tau,s} |< R/2),
$$
$$
P_{\cF_{t}}(\sup_{s\leq R^{2}}|m_{\tau,s} |< R/2)
\geq 2\bar \xi(d,\delta)>0.
$$
Also note that by Remark \ref{remark 10.14.1}
$$
P_{\cF_{\tau}}( \theta_{\tau}\tau'_{R}  < R^{2} )\leq
P_{\cF_{\tau}}\big(\int_{0}^{\theta_{\tau}\tau_{ R }}|b_{\tau+s}| \,ds
\geq R/2\big)
+
P_{\cF_{\tau}}(\sup_{s\leq R^{2}}|m_{\tau,s} |\geq R/2)
$$
$$
\leq 2\bar b_{R}+1-P_{\cF_{\tau}}(\sup_{s\leq R^{2}}|m_{\tau,s} |< R/2)\leq 2\bar b_{R}+1-2\bar\xi 
$$
and we get the second relation  in \eqref{8.2.2}  for 
 $2\bar b_{R}\leq 
\bar\xi$.

To prove the first relation take $\zeta$ such that
$\zeta(z)=\eta(z/R)$, where $\eta(z)=1$
 for $|z|\leq 2 $ and take
$T=\delta^{2}R^{2} $, in which case $u(0,x)\leq u(0,0)<1$
and $u(0,0)$ depends only on $\delta$ (and $\eta$).
Also define $\mu$ as the first time $\eta_{\tau,s}$
reaches $T$, which is certainly less than or
equal to $R^{2}$. Now observe that
$u(\eta_{\tau,\mu},|y+m_{\tau,\mu}|)=u(T,|y+m_{\tau,\mu}|)=\zeta
(|y+m_{\tau,\mu}|)$. It follows that
 $$
 P_{\cF_{\tau}}(\sup_{s\leq R^{2}}|y+m_{\tau,s}|< 2R)
\leq P_{\cF_{t}}( |y+m_{\tau,\mu}|< 2R)
$$
$$
\leq E_{\cF_{t}}u(\eta_{\tau,\mu},|y+m_{\tau,\mu}|)\leq u(0,|y|)\leq u(0,0).
$$
Hence,
$$
P_{\cF_{\tau}}( \theta_{\tau}\tau'_{ R} (y)  < R^{2} )\geq
 P_{\cF_{\tau}}\big(\int_{\tau}^{\tau+\theta_{\tau}\tau_{ R} (y) }|b_{s}| \,ds
\leq R/2,\sup_{s\leq R^{2}}|y+m_{\tau,s}|\geq 2R \big)
$$
$$
\geq 1-P_{\cF_{\tau}}\big(\int_{\tau}^{\tau+\theta_{\tau}\tau_{ R} (y)}|b_{s}| \,dt
\geq R/2\big)-P_{\cF_{\tau}}(\sup_{s\leq R^{2}}|y+m_{\tau,s}|\leq 2R)
$$
and it is clear how to adjust \eqref{12.18.3}
to get both inequalities in \eqref{8.2.2} with perhaps
$\bar\xi$ different from the above one.

To prove \eqref{1.3.1}   observe that
in light of \eqref{8.2.2}  for any $i=0,1,2...$
$$
P_{\cF_{\tau+iR^{2}}}(\max_{s\leq R^{2}}|x_{s+ \tau+iR^{2}}
-x_{ \tau+iR^{2}}+\xi_{i}|<R^{2})\leq 1-\bar\xi,
$$
where $\xi_{i}=x_{_{\tau+iR^{2}}}-x_{\tau}-y $. In other words,
$$
P_{\cF_{\tau+iR^{2}}}(\max_{s\leq R^{2}}|x_{s+ \tau+iR^{2}}
-x_{\tau}-y|<R^{2})\leq 1-\bar\xi.
$$
Now \eqref{1.3.1} follows since its left hand side is
the conditional expectation given $\cF_{t}$
of the product of the above probabilities over $i=0,..,n-1$.

To prove \eqref{1.3.3} note that  
$$
I=\sum_{n=1}^{\infty}E_{\cF_{\tau}}I_{\theta_{\tau}\tau'_{  R} (y)
>(n-1)R^{2} }E\Big\{
\int_{(n-1)R^{2} }^{(nR^{2})\wedge\theta_{\tau}\tau'_{ R} (y)} |
b_{\tau+s}|\,ds\mid \cF_{\tau+(n-1)R^{2}  }\Big\}
$$
$$
\leq \bar b_{R}R\sum_{n=1}^{\infty}P_{\cF_{\tau}}( \theta_{\tau}\tau'_{  R} (y)
>(n-1)R^{2} )
 \leq \bar b_{R}R \sum_{n=1}^{\infty}
(1-\bar\xi)^{n-1}.
$$
This yields \eqref{1.3.3}.  

To prove \eqref{1.2.1} we may assume that $|x|>R/16$ and
using assertion (iii) of Lemma \ref{lemma 1.2.1}
 conclude that
$$
du(| x_{\tau+s}-x_{\tau}-y|)\geq b^{i}_{\tau+s}D_{i}u(|x_{\tau+s}-x_{\tau}-y|)\,ds+dM_{s},
$$
where $M_{s}$ is a local martingale
before $x_{\tau+s}-x_{\tau}-y$ hits the origin. 
For our $x_{\cdot}$, on the time interval, which we denote
 $(0,\nu)$, when $ x_{\tau+s}-x_{\tau}-y
\in B_{R}\setminus \bar B_{R/16}$ we have
$|D u(|x_{\tau+s}-x_{\tau}-y|) \leq N(d,\alpha)R^{-\alpha-1}$.  
Furthermore, at starting point $u(-y)\geq (9R/16)^{-\alpha}$.
Consequently and by \eqref{1.3.3}
$$
(9R/16)^{-\alpha}\leq NR^{-\alpha-1}
E_{\cF_{\tau}}\int_{0}^{\theta_{\tau}\tau'_{  R} (y) }|b_{\tau+s}|\,ds
+P_{\cF_{\tau}}\big(\nu=\theta_{\tau}\tau'_{  R} (y)  \big)R^{-\alpha}
$$  
$$
+
P_{\cF_{\tau}}\big(\nu=\theta_{\tau}\gamma'_{ R/16}(x)  \big)(R/16)^{-\alpha},
$$
$$
(16/9)^{\alpha}\leq N \bar b_{R}
+1+(16^{\alpha}-1)
P_{\cF_{\tau}}\big(\theta_{\tau}\tau'_{  R} (y) >\theta_{\tau}\gamma'_{ R/16}(y)  \big) .
$$
It follows easily that \eqref{1.2.1} holds
with $\bar\xi$ perhaps different from the above ones,
once a relation like \eqref{12.18.3} holds.
The continuity of $m_{b}$ in \eqref{12.18.3}
and of $\bar\xi $
with respect to $\delta$ is established by inspecting
the above proof.
The theorem is proved.  \qed

 We remind the reader that from this point on
throughout the paper
we suppose that Assumption   $(\bar  b_{\rho_{b}})$  \ref{assumption 8.19.2} is satisfied. 
Theorem \ref{theorem 8.2.1} will allow us to show that the behavior of $x_{\tau+s}-x_{\tau}$
 before it exits from $B_{\rho_{b}}$ is almost the same as if there were no drift.

In light of \eqref{3.7.2} and \eqref{1.3.3} 
estimate \eqref{3.7.1} implies the following.

\begin{corollary}
                                     \label{corollary 3.7.1}
For any Borel $g\geq0$ and $\cF_{\tau}$-measurable $\bR^{d}$-valued $y$
we have
\begin{equation}
                                          \label{9.29.20}
 E_{\cF_{\tau}}\int_{0}^{\theta_{\tau}\tau'_{ R}(y)  }  
g( x_{\tau+s})\,ds\leq 
 N(d,  \delta)(1+\bar b_{R})  
 R 
\|g\|_{L_{d  }(\bR^{d})}.
\end{equation}
\end{corollary}

From \eqref{1.3.1} we immediately obtain the following
  \begin{corollary}
                                  \label{corollary 2.3.1}
Let   $R\in(0,\rho_{b}] $. Then there
exists a constant $N$, depending only on
$\bar\xi $, such that, for any $\cF_{\tau}$-measurable $\bR^{d}$-valued $y$,
$T \geq 0$,
$$
P_{\cF_{\tau}}(\theta_{\tau}\tau'_{R}(y)   >T)\leq Ne^{-T/(NR^{2})}.
$$
\end{corollary}

Next corollary serve as a lemma for Theorem \ref{theorem 8.20.1}.

\begin{corollary}
                                          \label{corollary 8.18.1}
For  $\mu \in[0,1]$,   stopping time $\tau$,  and any $R\in(0,\rho_{b}]$  we have
\begin{equation}
                                          \label{8.18.3}
I_{\tau<\infty}E_{\cF_{\tau}}e^{-\mu R^{-2} \theta_{\tau}\tau_{ R }} \leq e^{- \mu\bar \xi/2}.
\end{equation}
\end{corollary} 

 Indeed, for $\gamma:=\theta_{\tau}\tau_{ R} $
on the set $\{\tau<\infty\}$
the derivative with respect to $\mu$ of the left-hand
side of \eqref{8.18.3} is
$$
-R^{-2} E_{\cF_{\tau}} \gamma e^{-R^{-2}\mu \gamma} \leq -
 e^{-\mu } 
P_{\cF_{\tau}}( \gamma\geq R^{2})\leq-e^{-\mu}\bar \xi,
$$
where the last inequality follows from \eqref{8.2.2}.
By integrating we find
$$
E_{\cF_{t}}e^{-\mu R^{-2} \gamma} -1\leq  
(e^{-\mu   }-1)
\bar \xi,
$$
which after using
$$
e^{-\mu  }-1\leq- \mu /2,\quad
1-\mu
\bar \xi/2\leq e^{-\mu 
\bar \xi/2}
$$
leads to \eqref{8.18.3}.
 
\begin{theorem}
                                     \label{theorem 8.20.1}
For any $\lambda\geq \rho_{b}^{-2} $, stopping time $\tau$, and $R\in(0,\infty)$ on the set $\{\tau<\infty\}$ we have
\begin{equation}
                           \label{8.20.1}
E_{\cF_{\tau}}e^{-\lambda  \theta_{\tau}\tau_{ R} }\leq
e^{\bar\xi/2}e^{- \sqrt{\lambda}  R
\bar\xi/2} .
\end{equation}
   
In particular, if $s\leq R \rho_{b}
\bar \xi/4 $ we have
\begin{equation}
                            \label{10.2.2}
P_{\cF_{\tau}}(  \theta_{\tau}\tau_{ R} \leq  s )\leq 
 e^{\bar\xi/2}\exp\Big(-\frac{{\bar \xi}^{2}R^{2}}{16 s}\Big).
\end{equation}
\end{theorem}

Proof. In case $\tau<\infty$ take an integer $n\geq 1$ and introduce $\tau^{k}$, $k=1,...,n$,
as the first exit time of $(\tau+s,x_{\tau+s} )$, $s\geq0$,
from $C_{R/n}( \tau+\tau^{k-1 } ,x_{\tau+\tau^{k-1}} )$ after $\tau^{k-1}$, provided that   $\tau^{k-1}<\infty$
($\tau^{0}:=0$). 
If
 $
\lambda\leq n^{2}/R^{2}$ and $ R/n\leq \rho_{b}$ 
then by \eqref{8.18.3} with $\mu=(R/n)^{2}
\lambda$ we have
$$
I_{\tau^{k-1}<\infty}E_{\cF_{\tau^{k-1}}} e^{-\lambda( \tau^{k} - \tau^{k-1} )} \leq e^{- (R/n)^{2}
\lambda \bar \xi/2}.
$$
Hence,
\begin{equation}
                                           \label{12.27.1}
E_{\cF_{\tau}}e^{-\lambda \theta_{t}\tau_{ R} }\leq E_{\cF_{\tau}}\prod_{k=1}^{n}
e^{-\lambda( \tau^{k} - \tau^{k-1}  )}\leq  
 e^{- R^{2}n^{-1}
\lambda \bar \xi/2}.
\end{equation}
By taking
$n=\lceil R\sqrt\lambda\rceil$ and observing that that
$ R/n\leq  \rho_{b}$ and
$R^{2}n^{-1}
\lambda\geq R\sqrt\lambda-1$, we come to \eqref{8.20.1}.  

To prove \eqref{10.2.2} again consider the case that $\tau<\infty$ and
note that  for $\lambda\geq\rho_{b}^{-2}$ we have
$$
P_{\cF_{\tau}}( \theta_{\tau}\tau_{ R}   \leq  s  )=P_{\cF_{\tau}}\big(
 \exp(-\lambda 
 \theta_{\tau}\tau_{ R}  ) \geq \exp(-\lambda s)\big)
\leq  
 \exp(\bar \xi/2+ \lambda s-\sqrt{\lambda}R \bar \xi/2).
$$
For $\sqrt{\lambda}=R\bar \xi/(4s)$
 we get \eqref{10.2.2} provided that $R\bar\xi/(4s)\geq\rho_{b}^{-1}$.
The theorem is proved. \qed

\begin{corollary}
                                       \label{corollary 7.29.1}
Let   $\lambda>0$, $R\in (0,\rho_{b}]$, $t\geq0$.
Then there are constants $N=N(\bar\xi),\nu=\nu(\bar\xi)>0$ such that on the set $\{\tau<\infty\}$ we have
\begin{equation}
                                                   \label{8.21.1}
N E_{\cF_{\tau}}\int_{0}^{\theta_{\tau}\tau_{ R}}
e^{-\lambda t} \,dt  \geq \lambda^{-1}
(1-e^{-\lambda\nu R^{2}}) .
\end{equation}
In particular (as $\lambda\downarrow0$),
$NE_{\cF_{\tau}}\theta_{\tau}\tau_{ R}\geq \nu R^{2}$.

\end{corollary}

Indeed, for any $\nu\leq   \bar \xi/4$
we have
$$
E_{\cF_{\tau}}\int_{0}^{\theta_{\tau}\tau_{ R}}
e^{-\lambda t} \,dt=\lambda^{-1}
E_{\cF_{\tau}}(1-e^{-\lambda \theta_{\tau}\tau_{C_R}})
$$
$$
\geq\lambda^{-1}
 E_{\cF_{\tau}}I_{ \theta_{\tau}\tau_{ R} >\nu R^{2}} (1-e^{-\lambda\nu R^{2}})
=\lambda^{-1} P_{\cF_{\tau}}( \theta_{\tau}\tau_{ R} >\nu R^{2})(1-e^{-
\lambda\nu R^{2}})
$$
$$
\geq \lambda^{-1}\Big(1-e^{\bar\xi/2} \exp
\Big(-\frac{\bar\xi^{2} }{16\nu }\Big)\Big)
(1-e^{-\lambda\nu R^{2}}),
$$
which yields \eqref{8.21.1} for an appropriate 
small 
$\nu =\nu( \bar\xi )>0$.

This result will be used in proving a higher summability
of the Green's functions of $x_{\cdot}$.
The next one is aimed at proving the precompactness
of distributions of various processes like $x_{\cdot}$.

\begin{corollary}
                           \label{corollary 10.26.1}
For any $n>0$ and
  $t\geq 0$ on the set $\{\tau<\infty\}$ we have
\begin{equation}
                                  \label{10.28.2}
E_{\cF_{\tau}}\sup_{r\in[0,t]}|x_{\tau+r}-x_{\tau}|^{ n}
\leq N(  t ^{ n/2}+ t ^{ n}),
\end{equation}
where $N=N(n, \rho_{b},\bar\xi)$.
\end{corollary}

Indeed,  
for   $t\leq \mu \rho_{b}\bar\xi/4$
on the set $\{\tau<\infty\}$
we have
$$
P_{\cF_{\tau}}(\sup_{r\leq t }|x_{\tau+r}-x_{\tau}|\geq \mu)
\leq P_{\cF_{\tau}}(\theta_{\tau}\tau_{ \mu}\leq  t)
\leq e^{\bar\xi/2}\exp\Big(-\frac{  \mu^{2}\bar\xi^{2}}
{16 t}\Big).
$$
Consequently,
$$
E_{\cF_{\tau}}\sup_{r\leq t}|x_{\tau+r}-x_{\tau} |^{ n}
=n\int_{0}^{\infty}\mu^{n-1}
P_{\cF_{\tau}}(\sup_{r\leq t }|x_{\tau+r}-x_{\tau}|\geq \mu)\,d\mu
$$
$$
\leq n\int_{0}^{4(s-t)/(\rho_{b}\bar\xi)}\mu^{n-1}\,d\mu
+ne^{\bar\xi/2}\int_{0}^{\infty}\mu^{n-1}
\exp\Big(-\frac{  \mu^{2}\bar\xi^{2}}
{16 t}\Big)
\,d\mu,
$$
and   the result follows. 

A few more general results are related to going through
a long ``sausage".
\begin{theorem} 
                                        \label{theorem 1.24.1}
Let $R\in(0,\rho_{b}]$, $\tau$ be a   stopping time, $\cF_{\tau}$-measurable $y_{\tau}\in\bR^{d}$ be such that $16|x_{\tau}-y_{\tau}|\geq 3R$ on the set $\{\tau<\infty\}$.  On the same set
for $r>0$ denote by $S_{r}(x_{\tau},y_{\tau})$   the open convex hull
of $B_{r}(x_{\tau})\cup B_{r}(y_{\tau})$. Then there exist
$T_{0},T_{1}$, depending only on $\bar\xi$,
such that $0<T_{0}<T_{1}<\infty$ and on the set $\{\tau<\infty\}$ the $P_{\cF_{\tau}}$-probability $\pi$
that $x_{\tau+s}$, $s\geq0$, will reach $\bar B_{R/16}(y_{\tau})$ before exiting
from $S_{R}(x_{\tau},y_{\tau})$ and this will happen
on the time interval $[nT_{0}R^{2},nT_{1}R^{2}]$
is greater than $\pi_{0}^{n}$, where
$$
n= \Big\lfloor \frac{16|x_{\tau}-y_{\tau}|+R}{4R}\Big\rfloor ,\quad \pi_{0}=\bar\xi/3.
$$

\end{theorem}

Proof.  We argue in case $\tau(\omega)<\infty$.
Introduce $\nu=\nu(x_{\tau},y_{\tau}) $ as the first time $ x_{\tau+s}$
reaches $\bar B_{R/16}(y_{\tau})$ and $\gamma=\gamma(x_{\tau},y_{\tau}) $ as the first time
it exits from $S_{R}(x_{\tau} ,y_{\tau})$. Owing to
$16|x_{\tau}-y_{\tau} |\geq 3R$, we have $n\geq1$ and we are going to use the induction
on $n$ with the induction hypothesis that,
for all $R\in(0,\rho_{b}]$,
$$
\Big\lfloor \frac{16|x_{\tau}-y_{\tau} |+R}{4R}\Big\rfloor=n
\Longrightarrow 
P_{\cF_{\tau}}(\gamma >\nu \in[nT_{0}R^{2},n
T_{1}R^{2}])\geq \pi^{n}_{0}.
$$

   If $n=1$, then $3R/16\leq |x_{\tau}-y_{\tau}|< 7R/16$ and by Theorem
\ref{theorem 8.2.1} (see \eqref{1.2.1}) we have $P_{\cF_{\tau}}(\theta_{\tau}\tau'_{ R } >\nu )\geq\bar\xi$.
Furthermore, in light of  Theorem
\ref{theorem 8.2.1}, there is $T_{1}=T_{1}(\bar\xi)$
such that $P_{\cF_{\tau}}(\theta_{\tau}\tau'_{ R } >T_{1}R^{2})\leq \bar\xi/3$.
Using \eqref{10.2.2} we also see that there is
$T_{0}=T_{0}(\bar\xi)<T_{1}$ such that
$P_{\cF_{\tau}}(\nu \leq T_{0}R^{2})\leq \bar\xi/3$.
It follows that $P_{\cF_{\tau}}(\gamma >\nu \in[T_{0}R^{2},
T_{1}R^{2}])\geq\bar\xi/3=\pi_{0}$.
This justifies the start of the induction.

Assuming that our hypothesis is true for some $n\geq 1$
 suppose that
$(n+2)R/4>|x_{\tau}-y_{\tau}|+R/16\geq  (n+1)R/4$. In that case, let 
$$
z_{\tau}=nR(x_{\tau}-y_{\tau})/(4|x_{\tau}-y_{\tau}|)
$$
and let $\nu'$ be the first time $ x_{\tau+s}$ reaches $\bar B_{R/16}(z_{\tau})$,
and let $\gamma'$ be the first time it exits
from $S_{R}(x_{\tau} ,z_{\tau})$. As is easy to see,
$$
P_{\cF_{\tau}}(\gamma >\nu\in[(n+1)T_{0}R^{2},(n+1)
T_{1}R^{2}])
$$
$$
\geq P_{\cF_{\tau}}\big(\gamma'>\nu'\in[T_{0}R^{2},T_{1}R^{2}],  
\gamma(x_{\nu'},y_{\tau}) >\nu(x_{\nu'},y_{\tau}) \in[nT_{0}R^{2},n
T_{1}R^{2}] \big)
$$
$$
=E_{\cF_{\tau}}I_{\gamma'>\nu'\in[T_{0}R^{2},T_{1}R^{2}]}
P_{\cF_{\nu'}} \big(\gamma(x_{\nu'},y_{\tau}) >\nu(x_{\nu'},y_{\tau})\in[nT_{0}R^{2},n
T_{1}R^{2}] \big).
$$
Observe that on the set $\nu'<\infty$ we have
$nR/4\leq |x_{\nu'}-y_{\tau}|+R/16<(n+1)R/4$, so that, by  our induction hypothesis, the last conditional probability
above is greater than $\pi_{0}^{n}$. Then   using the first part of the proof
we obtain our result for $n+1$ in place of $n$.
The theorem is proved. \qed

\begin{remark}
                                      \label{remark 1.25.1}
Notice that, for any fixed $x_{\tau},y_{\tau}$, the interval
$[nT_{0}R^{2},nT_{1}R^{2}]$ is as close to zero
as we wish if we choose $R$ small enough.
Then, of course, the corresponding probability will be
quite small but $>0$.

\end{remark} 

\begin{corollary}
                                    \label{corollary 10.17.1}
Let $R\in(0,\infty)$,  $\tau$ be a  stopping time
and let $y$ be $\cF_{\tau}$-measurable $\bR^{d}$-valued. Then
there is a constant $N=N(d,\delta,\rho_{b},R)$ such that
on $\{\tau<\infty\}$ for all $T\in(0,\infty)$ we have
\begin{equation}
                                            \label{10.17.3}
P_{\cF_{\tau}}(\max_{t\leq T}|x_{\tau+t}-y|<R)
\leq Ne^{-T/N}.
\end{equation}

\end{corollary}

Indeed, on the set where $\tau<\infty$ and $|x_{\tau}-y|
\geq R$ estimate \eqref{10.17.3} is obvious. So we may
concentrate on $y\in B_{R}(x_{\tau})$ Then
by using 
Theorem \ref{theorem 1.24.1} we see that, given that
$\tau<\infty$ with    
 probability not less than some $\beta>0$,
  depending only on $d,\delta$, and $\rho_{b}$,
the process   $ x_{\tau+t}$
 will reach $\bar B_{\rho_{b}/16}(x_{\tau} +\rho_{b}e_{1}/4)$,
where $e_{1}$ is the first basis vector, before time $T_{1}=T_{1}
(d,\delta)$. Therefore, its first coordinate will increase
by at least $3\rho_{b}/16$.   Repeating this argument
and taking into account that $R<\infty$, we see that
with probability
$\pi>0$ depending only on $d,\delta$, $\rho_{b}$, and $R $,
the process $ x_{\tau+t} $ will leave $ B_{R}(y)$ before time
$S$, where $S$ depends only on $d,\delta$, $\rho_{b}$, and $R $,  
that is
$$
P_{\cF_{\tau}}(\max_{t\leq S}|x_{\tau+t}-y|<R)\leq 1-\pi.
$$
Iterating this inequality, which is also true if
$|x_{\tau}-y|
\geq R$, we obtain for $n=1,2,...$  that $
P_{\cF_{\tau}}(\max_{t\leq nS}|x_{\tau+t}-y|<R)\leq (1-\pi)^{n}$ and this yields the result.

The following complements Corollaries \ref{corollary 2.3.1}
and \ref{corollary 10.17.1}.
 
\begin{corollary}
                                     \label{corollary 1.25.1}
Let   $\kappa\in[0,1)$, $R\in(0,\rho_{b}]$, and $|x|\leq\kappa R$ and let $\tau$ be a stopping time.
Then for any $T>0$ on the set $\{\tau<\infty\}$
\begin{equation}
                                        \label{1.25.2}
NP_{\cF_{\tau}}(\theta_{\tau}\tau' _{ R}(x)  > T)\geq e^{-\nu T/[(1-\kappa)R]^{2}},
\end{equation}
where $N$ and $\nu>0$ depend only on $\bar\xi$.
\end{corollary}

Indeed, passing from $B_{R}(x)$ to $B_{(1-\kappa)R} \subset
B_{R}(x)$
shows that we may assume that $x=0$ and $\kappa=0$.
In that case, consider meandering of $x_{\tau+s}$
between $\bar B_{R/16}$ and $\partial B_{R/16}(y)$, where
$|y|=R/4$, without exiting from $B_{R}$. As is easy to deduce
from Theorem \ref{theorem 1.24.1},
given that the $n$th loop happened, with probability
$\pi^{4}_{0}$ the next loop will occur and take at least
$4R^{2}T_{0}$ of time. Thus  the $n$th loop
will happen and will take   at least $4nR^{2}T_{0}$ of time
with probability at least $\pi_{0}^{4n}$.
It follows that, for any $n$,
$$P_{\cF_{\tau}}(\theta_{\tau}\tau' _{ R  }  \geq 4nR^{2}T_{0}))\geq \pi_{0}^{4n},
$$
and this yields \eqref{1.25.2} for $x=0$ and $\kappa=0$.

\def\lpq{$L_{{p,q}}\,$}

\mysection{First estimates of potentials in $L_{p,q}$
}

Some arguments in the future have to be repeated
  twice in slightly different
situations. In order to avoid this
we consider the following setting.

Take some   $p,q\in(1,\infty]$, $\alpha \in[0,1]$
and fix $ \kappa(\bar\xi) >1 $ such that
the right-hand side of \eqref{8.20.1} equals $1/2$ when $R=(\kappa(\bar\xi)-1)\lambda^{-1/2}$.

 \begin{assumption} 
                            \label{assumption 8.24.1}
There exists a constant $\ell$ such that
 for
any $\lambda\geq \kappa^{2}(\bar\xi)\rho_{b}^{-2}$,   stopping time $\tau$, and Borel $f\geq0$ given on $\bR^{d+1}$, it holds that
\begin{equation}
                              \label{8.24.1}
 E_{\cF_{\tau}}\int_{0}^{\theta_{\tau}\tau_{\rho(\lambda)}}e^{-\lambda s}f(\tau+s,x_{\tau+s})\,ds\leq \ell\lambda ^{-\alpha}\|f\|_{L_{p,q}},
\end{equation}
 where $\rho(\lambda)=\kappa(\xi)\lambda^{-1/2} $.
\end{assumption}

 \begin{remark}
                        \label{remark 8.24.1}  
As it follows from Theorem \ref{theorem 9.27.1} with
$ \gamma=\tau+\theta_{\tau}\tau_{\rho(\lambda)}$, Assumption \ref{assumption 8.24.1}
is satisfied, for instance, if $p=q=d+1$     
   with $\ell \lambda^{-\alpha}$ replaced by
$$
 N(d,\delta)(\lambda^{-1}+\hat b^{2}_{\rho(\lambda)})^{d/(2d+2)}
\leq N(d,\delta)(1+\bar b^{2}_{\rho_{b}})^{d/(2d+2)}\lambda^{-d/(2d+2)}
$$
$$
\leq N(d,\delta)2^{d/(2d+2)}\lambda^{-d/(2d+2)}\leq N(d,\delta)\lambda^{-d/(2d+2)}.
$$

Another case is when we take $p=d$, $q=\infty$. 
This time Assumption \ref{assumption 8.24.1}
satisfied with $\ell \lambda^{-\alpha}$ replaced by
$$
 N(d,\delta)(\lambda^{-1}+\hat b^{2}_{\rho(\lambda)})^{1/2}
\leq N(d,\delta)(1+\bar b^{2}_{\rho_{b}})^{1/2}\lambda^{-1/2}
$$
$$
\leq N(d,\delta)2^{1/2}\lambda^{-1/2}\leq N(d,\delta)\lambda^{-1/2}.
$$

 Below we will see that there is a wider
range of parameters $p,q$.
 \end{remark}

\begin{remark}
                        \label{remark 10.15.1}
Since the $L_{p,q}$-norm is translation invariant,
one gets an equivalent assumption if 
$(\tau+s,x_{\tau+s})$ is replaced with  
$(\nu+s ,x_{\tau+s}+y)$ for any $\cF_{\tau}$-measurable
real-valued
$\nu$ and $\bR^{d}$-valued $y$. In short, in Assumption \ref{assumption 8.24.1} one can replace \eqref{8.24.1} with
\begin{equation}
                              \label{10.15.1}
 E_{\cF_{\tau}}\int_{0}^{\theta_{\tau}\tau_{\rho(\lambda)}}e^{-\lambda s}f( \nu+s,x_{\tau+s}+y)\,ds\leq \ell\lambda ^{-\alpha}\|f\|_{L_{p,q}}.
\end{equation}

\end{remark}

\begin{lemma}
                                     \label{lemma 8.22.1}

 Under Assumption \ref{assumption 8.24.1} 
for any  stopping time $\tau$, $s_{0}\geq0$, $y_{0}\in \bR^{d}$, $\lambda\geq \kappa^{2}(\bar\xi)\rho_{b}^{-2} $,
and Borel nonnegative $f$ vanishing outside 
$C_{\lambda^{-1/2}}(s_{0},y_{0})$
we have on the set $\{\tau<\infty\}$ that
\begin{equation}
                                             \label{8.22.10}
 E_{\cF_{\tau}}\int_{0}^{\infty}e^{- \lambda s}
f(s,x_{\tau+s}-x_{\tau})  \,ds\leq 
N(\bar\xi)\ell\lambda^{-\alpha } 
\Phi_{\lambda}(s_{0},y_{0})\|  f\|_{L_{p ,q }},
\end{equation}
where   $\Phi_{\lambda}
(t,x)=e^{-\sqrt{\lambda}  (\sqrt t+|x|)\bar\xi/4}$.
 
\end{lemma}

Proof. We assume that the event $\{\tau<\infty\}$
occurred. Introduce
  $\tau^{0}$ as the first time $(s,x_{\tau+s}-x_{\tau})$, $s\geq0$,
hits $\bar C_{\lambda^{-1/2}}(s_{0},y_{0})$ and set $\gamma^{0}$
as the first time after $\tau^{0}$ this process
 exits from $C_{ 
\kappa(\xi)\lambda^{-1/2}  }(s_{0},y_{0})$ after $\tau^{0}$.
We define recursively $\tau^{k}$, $k=1,2,...$, as the first
time after $\gamma^{k-1}$ the process
$(s-\tau,x_{s}-x_{\tau})$
hits $\bar C_{\lambda^{-1/2}}(s_{0},y_{0})$ and $\gamma^{k}$ as 
the first time after $\tau^{k}$ this process
 exits from 
$C_{ \kappa(\xi)\lambda^{-1/2} }(s_{0},y_{0})$.
As is easy to see $\tau+\tau^{k}$ and $\tau+\gamma^{k}$
are stopping times and they  are either
infinite or lie between $\tau+s_{0}$ and $\tau+s_{0}+\kappa^{2}(\xi)\lambda^{-1}$.

The left-hand side of
\eqref{8.22.10} on $\{\tau<\infty\}$ equals
\begin{equation}
                                                \label{8.22.1}
 E_{\cF_{\tau}}
\sum_{k=0}^{\infty} e^{- \lambda \tau^{k}  }I_{k},
\end{equation}
where
$$
I_{k}=I_{\tau^{k}<\infty }
E_{\cF_{\tau+\tau^{k}}} \int_{\tau^{k} }^{ \gamma^{k} 
 }
e^{-\lambda  (s-\tau^{k}) }
f( s,x_{\tau +s}-x_{\tau}) \,ds  
$$
$$
=I_{\tau^{k}<\infty }E_{\cF_{\tau+\tau^{k}}} \int_{0 }^{\gamma^{k}-\tau^{k}
 }
e^{-\lambda r }
f( r+\tau^{k} ,x_{\tau+\tau^{k}+r}-x_{\tau}) \,dr .
$$
By Remark \ref{remark 10.15.1} the last expression
is less than $\ell\lambda^{-\alpha}\|f\|_{L_{  p , q }}$,
that is $I_{k}\leq \ell\lambda^{-\alpha}\|f\|_{L_{  p , q }}$.

Next,   observe that, if $\sqrt{s_{0}}>|y_{0}| $, then  $\tau^{0} $
is bigger than the first exit time of $( s,x_{\tau+s}-x_{\tau})$
from $C_{\sqrt{s_{0}}}$, that is $\tau^{0} \geq
\theta_{\tau}\tau_{ \sqrt{s_{0} }}$ and
by Theorem \ref{theorem 8.20.1} 
$$
E_{\cF_{\tau}}e^{- 
\lambda  \tau^{0}  }\leq N(\bar \xi)e^{-\sqrt{\lambda} 
\sqrt{s_{0}}\bar\xi/2}.
$$
In case $\sqrt{s_{0}}\leq |y_{0}| $ and $|y_{0}|> \lambda^{-1/2}$
our $\tau^{0}$
is bigger than $\theta_{\tau}\tau_{ |y_{0}|-\lambda^{-1/2}} $ and
$$
E_{\cF_{\tau}}e^{- 
\lambda  \tau^{0}  }\leq Ne^{-\sqrt{\lambda}(|y_{0}|-\lambda^{-1/2})\bar\xi/2}.
$$
The last estimate (with   $N=1$) also holds if 
$|y_{0}|\leq\lambda^{-1/2}$, so that
in case $\sqrt{s_{0}}\leq |y_{0}| $
$$
E_{\cF_{\tau}}e^{- 
\lambda  \tau^{0}  }
\leq Ne^{-\sqrt{\lambda} |y_{0}|\bar\xi/2}
$$
and we conclude that in all cases
$$
E_{\cF_{\tau}}e^{- 
\lambda  \tau^{0} } \leq 
N(\bar\xi)e^{-\sqrt{\lambda} (\sqrt{t_{0}}+|y_{0}|)\bar\xi/4}.
$$

Furthermore,  by the choice of $\kappa(\bar\xi)$ and
 Theorem \ref{theorem 8.20.1}
$$
E_{\cF_{\tau+\gamma^{k-1}}}e^{-\lambda ( \tau^{k}
- \gamma^{k-1} )} \leq \frac{1}{2},
$$
$$
E_{\cF_{\tau}}e^{-\lambda  \tau^{k }}=
E_{\cF_{\tau}}e^{-\lambda  \gamma^{k-1}  }E_{\cF_{\tau+\gamma^{k-1}}}
 e^{-\lambda ( \tau^{k} 
- \gamma^{k-1} )}  
\leq \frac{1}{2}E_{\cF_{\tau}}e^{- \lambda  \gamma^{k-1}  },
$$
so that ($\gamma^{k-1}\geq \tau^{k-1}$)
$$
E_{\cF_{\tau}}e^{-\lambda   \tau^{k} }\leq\frac{1}{2}
 E_{\cF_{\tau}}e^{- \lambda  \tau^{k-1} } ,\quad
 E_{\cF_{\tau}}e^{- \lambda  \tau^{k} } \leq 2^{-k}
E_{\cF_{\tau}}e^{- \lambda  \tau^{0}  }.
$$

Recalling \eqref{8.22.1} we see that the left-hand side 
of \eqref{8.22.10}
is indeed dominated by
$
N\ell\lambda^{-\alpha}\Phi_{\lambda}(s_{0},y_{0})\|f\|_{L_{p ,q }}
$  
and the lemma is proved. \qed

 The following theorem shows that the time spent
by $(s,x_{s})$ in cylinders $C_{1}(0,x)$
decays very fast as $|x|\to\infty$. Define
$$
\bR^{d+1}_{+}=[0,\infty)\times\bR^{d}.
$$
\begin{theorem}
                               \label{theorem 8.22.1}
 Under Assumption \ref{assumption 8.24.1}
  for any  stopping time $\tau$, any $\lambda\geq \kappa^{2}(\bar\xi)\rho_{b}^{-2}$, and   Borel nonnegative $f$  
on the set $\{\tau<\infty\}$
we have 
\begin{equation}
                                     \label{8.22.4}
I:= E_{\cF_{\tau}}\int_{0}^{\infty}e^{- \lambda s}
f(s,x_{\tau+s}-x_{\tau}) \,ds\leq 
N (\bar\xi)\ell\lambda^{-\alpha }\|\Psi_{\lambda}  f\|_{L_{p ,q }(\bR^{d+1}_{+})},
\end{equation}
where $\Psi _{\lambda}(t,x)=\exp(- \sqrt{\lambda} 
(|x|+ \sqrt t)\bar\xi/16)$.
\end{theorem}

Proof.  
Set $\zeta(t,x)=\lambda^{ (d+2)/2 }\eta(\lambda t,\sqrt{\lambda} x)$, where $\eta$ has unit integral
and is proportional to the indicator of $C_{1}$,  and for $(t,x),(r,y)\in\bR^{d+1}$
set
$$
f_{r,y}(t,x)=f(t,x)\zeta(t-r,x-y).
$$
Clearly, due to Lemma \ref{lemma 8.22.1},
$$
I
=\int_{\bR}\int_{\bR^{d }}E_{\cF_{\tau}}\int_{0}^{\infty}
e^{-\lambda s}
f_{r,y}(s ,x_{\tau+s}-x_{\tau}) \,ds\,dydr
$$
$$
\leq N(\bar\xi)\ell\lambda^{-\alpha } \int_{-1/\lambda}^{\infty}\int_{\bR^{d }}
\Phi_{\lambda}(r_{+},y)\|  f_{r,y}\|_{L_{p ,q }}\,dydr.
$$
Further estimates of the  factor of $N(\bar\xi)\ell\lambda^{-\alpha }$ are done exactly as in the proof
of Theorem 2.9 of \cite{2022_2}. The theorem is proved. \qed

If  $f(t,x)=f(x)$ we come to the following.
\begin{corollary}
                           \label{corollary 10.4.1}
Under Assumption \ref{assumption 8.24.1}
for any $\lambda\geq \kappa^{2}(\bar\xi)\rho_{b}^{-2}$ and   Borel nonnegative $f(x)$  
we have on $\{\tau<\infty\}$
\begin{equation}
                                          \label{10.4.2}
 E_{\cF_{\tau}}\int_{0}^{\infty}e^{- \lambda s}
f( x_{\tau+s}-x_{\tau}) \,dt\leq 
N(\bar\xi)\ell\lambda^{-\alpha-1/q} \|\Psi_{\lambda}  f\|_{L_{p}(\bR^{d})},
\end{equation}
where $\Psi _{\lambda}( x)=\exp(- \sqrt{\lambda}  |x|\bar\xi/16)$.
\end{corollary}   
 
Next result  is dealing with the exit times of the process
$x_{s}$ rather than $(s,x_{s})$. We will need it
while showing an improved integrability of Green's functions. Assumption \ref{assumption 8.24.1}
is no longer needed.

Estimate \eqref{9.29.5} below in case $b$ is bounded
was the starting point for the theory
of {\em time homogeneous\/} controlled diffusion processes
about fifty years ago.
 \begin{lemma}
                                      \label{lemma 9.29.1}
Let $p \in [d,\infty]$ and let $\tau$ be a stopping time. Then for any Borel 
nonnegative $f(x)$, $R\leq \rho_{b}$, 
and $x\in\bR^{d}$  
\begin{equation}
                                       \label{9.29.5}
 E_{\cF_{\tau}}\int_{0}^{\theta_{\tau}\tau' _{ R} }   
f( x_{\tau+s})\,ds\leq N(\delta,d, \bar\xi)
 R 
 ^{2-d/  p  } 
\|f\|_{L_{p  }(\bR^{d})}.
\end{equation}

\end{lemma}

Proof. If $p=d$, the result follows from 
Theorem \ref{theorem 9.27.1}
and Remark \ref{remark 9.23.1}. Indeed,
in Theorem \ref{theorem 9.27.1} we have $A\leq N(R^{2}+B^{2})$
and, by definition, $B/R\leq \bar b_{R} \leq \bar b_{\rho_{b}}\leq m_{b}\leq1$.

In the general case we note as in Remark \ref{remark 10.15.1}
that 
$$
 E_{\cF_{\tau}}\int_{0}^{\theta_{\tau}\tau' _{ R} }   
f( x_{\tau+s}-x_{\tau})\,ds\leq N(\delta,d, \bar\xi)
 R 
\|f\|_{L_{d  }(\bR^{d})},
$$
where the last norm can be taken only over $B_{R}$ because
$x_{\tau+s}-x_{\tau}\in B_{R}$ before $\theta_{\tau}\tau' _{ R}$. After that we replace $R 
\|f\|_{L_{d  }(\bR^{d})}$ with $NR ^{2-d/p}
\|f\|_{L_{p }(\bR^{d})}$ by using H\"older's inequality
and then come back to \eqref{9.29.5} again by the argument
in Remark \ref{remark 10.15.1}. The lemma is proved. \qed

\mysection{Green's functions}  
                         \label{section 10.26.1}
In this section we  
take a stopping time $\tau$ and $\lambda\geq \kappa^{2}(\bar\xi)\rho_{b}^{-2}$, and introduce a measure on $\cF_{\tau}\otimes
\frB(\bR^{d+1})$ by
$$
G(\Gamma):= E I_{\tau<\infty}\int_{0}^{\infty}e^{- \lambda s}
I_{\Gamma}(\omega,s,x_{\tau+s}-x_{\tau}) \,ds.
$$
Take  $\cA\in\cF_{\tau}$ such that
$\cA\subset \{\tau<\infty \}$ and define
a measure on $\bR^{d+1}$ by
$$
G_{\cA}(\Lambda)= G(\cA\times\Lambda) ,
$$
Then for any Borel $f\geq 0$
on $\bR^{d+1}$ by
Remark \ref{remark 8.24.1}  and  
Theorem \ref{theorem 8.22.1} we have
$$
\int_{\bR^{d+1}}f\,G_{\cA}(dtdx)=
\int_{\Omega\times \bR^{d+1}}I_{A}f(t,x)\,G(d\omega dt dx)
$$
$$
=EI_{\cA }\int_{0}^{\infty}e^{- \lambda s}
f(s,x_{\tau+s}-x_{\tau}) \,ds
\leq 
N (d,\delta) \lambda^{-d/(2d+2) }\|\Psi_{\lambda}  f\|_{L_{d+1}(\bR^{d+1}_{+})}P(\cA ).
$$
This shows that the measure $G_{\cA}(\Lambda)$
has a density, which we denote by $G_{\cA}(t,x)$, and the following result holds true.
 
\begin{theorem}
                         \label{theorem 9.3.1}
The function $G_{\cA}$ is Borel measurable,
is such that $G_{\cA}(t,x)=0$ for $t\leq0$,
\begin{equation}
                                         \label{9.3.5}
 \|\Psi^{-1}_{\lambda} G_{\cA}\|_{L_{(d+1)/d } } 
\le N (d,\delta)\lambda^{ -d/(2d+2 )}P(\cA) ,
\end{equation}
and for any Borel nonnegative $f$ given on $\bR^{d+1}$   
we have
\begin{equation}
                                         \label{10.16.1}
EI_{\cA}\int_{0}^{\infty}e^{- \lambda s}
f(s,x_{\tau+s}-x_{\tau}) \,ds =\int_{\bR^{d+1}}f(t,x)G_{\cA}(t,x)\,dxdt.
\end{equation}
\end{theorem}

In light of \eqref{10.16.1} it is natural to call $G_{\cA}$
the Green's function of the process $(s,x_{\tau+s}-x_{\tau})$
on $\cA$.

 It turns out that,
actually,
$G_{\cA}$ is summable to a  power higher than $(d+1)/d$. The proof of this
is based on the parabolic version of
Gehring's lemma from \cite{GS_82} (also see Appendix 
in \cite{4}).

Introduce $\cC_{+}$ as the set of cylinders $C_{R}(t,x)$, $R>0$,  
$t\geq0$, $x\in\bR^{d}$. For $C=C_{R}(t,x)\in \cC_{+}$ let
 $2C=C_{2R}(t,x)$. If $C\in \cC_{+}$ and $C=C_{R}(t,x)$,
we call $R$ the radius of $C$.  Introduce
(for integrals in $\bR^{d+1}$ and $\bR^{d}$)
$$
\dashint_{\Gamma}f...=\big(\int_{\Gamma}1...)^{-1}
\int_{\Gamma}f...
$$

\begin{theorem}
                                  \label{theorem 9.3.2}
There exist  $d_{0}\in(1,d)$ and    $N $, depending only
on $\delta,d $,  
such that for any     $C\in \cC_{+}$ of radius 
$R\leq \kappa (\bar\xi)/(2\sqrt\lambda)$ and $p\geq d_{0}+1 $, we have
\begin{equation}
                                \label{10.14.01}
\| G_{ \cA}\|_{L_{p/(p-1)}(C)}\leq N R^{-(d+2) /p }
\| G_{\cA}\|_{L_{1}(2C )} , 
\end{equation}
which is equivalently rewritten as
$$
\Big(\dashint_{C}G^{p/(p-1)}_{\cA}\,dxdt\Big)^{(p-1)/p}
\leq N\dashint_{2C}G_{\cA}\,dxdt.
$$
 
\end{theorem}

Proof. We basically follow the idea in \cite{FS_84}.
Take $C\in \cC_{+}$ of radius $R\leq\kappa (\bar\xi)/(2\sqrt\lambda) $ and on the set $\{\tau<\infty\}$
define 
recursively  $\tau^{0}$ as the first  time
  when the process $(t ,x_{\tau+t}-x_{\tau})$, $t\geq0$,
hits $\bar C$, $\gamma^{0}$ as the first time after $\tau^{0}$
when this process leaves $2C$, $\tau^{n }$ as the first  time
after $\gamma^{n-1}$ when the process $(t ,x_{\tau+t}-x_{\tau})$
hits $\bar C$, $\gamma^{n }$ as the first time after $\tau^{n }$
when this process leaves $2C$.
 
Then for any nonnegative Borel $f$ vanishing outside $C$
with $\|f\|_{L_{d+1}(C)}=1$ 
we have
$$ 
I:=
\int_{C}f (t,x)  G_{\cA}(t,x)\,dxdt
$$  
$$
=\sum_{n=0}^{\infty}EI_{\cA}e^{- \lambda \tau^{n} } 
E_{\cF_{\tau+\tau^{n}}} \int_{0}^{\gamma^{n}-\tau^{n}}e^{-\lambda t }f(t+\tau^{n} ,x_{\tau+\tau^{n}+t}-x_{\tau})\,dt.
$$

Next, we use  
\eqref{9.29.2}  and what was said about the relation
of \eqref{8.24.1} to \eqref{10.15.1}  and take into account that $\bar b_{\kappa (\bar\xi)/\sqrt\lambda}\leq\bar b_{\rho_{b}}\leq 1$ 
 to see that on the set $\{\tau^{n}<\infty\}$ the conditional expectations
above are less than $N(d,\delta)R^{d/(d+1)}$. After that
we use   Corollary
\ref{corollary 7.29.1}   
to get that
$$
R^{2}I_{\tau^{n}<\infty}\leq N(\bar \xi) E_{\cF_{\tau+\tau^{n}}} \int_{\tau^{n}}^{\gamma^{n}}
e^{-\lambda(t-\tau^{n})} \,dt. 
$$
Then we obtain   
$$
I
\leq N  R ^{-(d+2)/(d+1)}
\sum_{n=1}^{\infty}EI_{\cA} e^{-\lambda \tau^{n} }
\int_{\tau^{n}}^{\gamma^{n}}
e^{-\lambda(t-\tau^{n})} \,dt
$$
$$
=N  R ^{-(d+2)/(d+1)}
\sum_{n=1}^{\infty}EI_{\cA}  
 \int_{\tau^{n} }^{\gamma^{n} }e^{-\lambda t } \,dt
\leq N  R ^{-(d+2)/(d+1)}
EI_{\cA} \int_{0}^{\infty}e^{-\lambda t}I_{2C}(t,x_{t}) \,dt
$$
$$
=N  R ^{-(d+2)/(d+1)}\int_{2C}G_{\cA}(t,x)\,dxdt.
$$

The arbitrariness of $f$ 
implies that
$$
\Big(\dashint_{C}  G_{\cA} ^{(d+1)/d}(t,x)\,dxdt
\Big)^{d/(d+1)}\leq N \dashint_{2C}  G_{\cA}(t,x)\,dxdt.
$$

Now the assertion of the theorem for $p=d_{0}+1$ follows directly from
the parabolic version of the famous Gehring's lemma
stated as Proposition 1.3  in \cite{GS_82}. 
For larger $p$ it suffices to use H\"older's inequality.
The theorem is proved. \qed

\begin{theorem}
                                           \label{theorem 2.3.1}
For any $p\geq d_{0}+1$  
\begin{equation}
                              \label{8.26.2}
\|G_{\cA}\|_{L_{p/(p-1)}(\bR^{d+1}_{+})}\leq N(\delta,d  )\lambda^{(d+2)/(2p)-1}P(\cA)
\end{equation}
In particular, for any Borel $f\geq0$ given on $\bR^{d+1}$ $($and $\lambda\geq \kappa^{2}(\bar\xi)\rho_{b}^{-2}$$)$ on $\{\tau<\infty\}$ we have
\begin{equation}
                               \label{8.26.5}
E_{\cF_{\tau}}\int_{0}^{\infty}e^{-\lambda t}f(t,x_{\tau+t})\,dt
\leq N(\delta,d )
\lambda^{(d+2)/(2d_{0}+2)-1}\|f\|_{L_{d_{0}+1}(\bR^{d+1})}.
\end{equation}
 
\end{theorem}

Proof.   Represent $\bR^{d+1}_{+}=[0,\infty)\times \bR$
as the union of countably many $C^{1},C^{2},...\subset \cC_{+}$
of radius $\kappa (\bar\xi)/(2\sqrt\lambda)$
so that each point in $\bR^{d+1}_{+}$ belongs  
to no more than $m(d)$  of the $2C^{i}$'s. Then 
$$
\|G_{\cA}\|_{L_{p/(p-1)}(\bR^{d+1}_{+})}
\leq \big\|\sum_{i}I_{C_{i}}G_{ \cA}\big\|_{L_{p/(p-1)}(\bR^{d+1}_{+})}
$$
$$
\leq\sum_{i} \| G_{\cA}\|_{L_{p/(p-1)}(C^{i})}
\leq N(\delta,d)\lambda^{(d+2)/(2p)}\sum_{i} \| G_{\cA}\|_{L_{1}(2C^{i})}
$$
$$
\leq N_{1}\lambda^{(d+2)/(2p)}
\| G_{\cA}\|_{L_{1}(\bR^{d+1}_{+})}= 
N_{1}\lambda^{(d+2)/(2p)-1}P(A).
$$
This proves \eqref{8.26.2} and the fact that
$$
EI_{\cA}\int_{0}^{\infty}e^{- \lambda t}
f(t,x_{\tau+t}-x_{\tau}) \,ds\leq 
N(\delta,d )
\lambda^{(d+2)/(2d_{0}+2)-1}\|f\|_{L_{d_{0}+1}(\bR^{d+1})}P(\cA).
$$
The arbitrariness of $\cA$ shows that \eqref{8.26.5}
holds with $f(t,x_{\tau+t}-x_{\tau})$ in place of 
$f(t,x_{\tau+t}-x_{\tau})$. One then eliminates 
$x_{\tau}$ as in Remark \ref{remark 10.15.1}.
The theorem is proved.

\begin{remark}
                             \label{remark 8.26.1}

If $\lambda\in(0,\kappa^{2}(\bar\xi)\rho_{b}^{-2})$, one can also give
an estimate of the left-hand side $J$ of \eqref{8.26.5} 
by taking nonnegative $f\in L_{p}(\bR^{d+1}_{+})$
and observing that, for $\lambda_{0}=\kappa^{2}(\bar\xi)\rho_{b}^{-2}$,
$$
J 
=\sum_{n=0}^{\infty}e^{-\lambda n}
E_{\cF_{\tau}} \int_{n}^{n+1}e^{-\lambda (t-n)}f(t,x_{\tau+t})\,dt
$$
$$
\leq \sum_{n=0}^{\infty}e^{\lambda_{0}-\lambda}e^{-\lambda n}
E_{\cF_{\tau}}\int_{n}^{n+1}e^{- \lambda_{0} (t-n)}f(t,x_{\tau+t})\,dt
$$
$$
=\sum_{n=0}^{\infty}e^{\lambda_{0}-\lambda}e^{-\lambda n}
E_{\cF_{\tau}}E_{\cF_{\tau+n}}\int_{0}^{ 1}e^{- \lambda_{0} t}f(n+t,x_{\tau+n+t})\,dt
$$
where  each conditional expectation in the sum is dominated
by 
$$
N\|fI_{[n,n+1)}\|_{L_{p }(\bR^{d+1}_{+})}
$$
 in light of \eqref{8.26.5} . Therefore
$$
J\leq N\sum_{n=0}^{\infty} e^{-\lambda n}\|fI_{[n,n+1)}\|_{L_{p }(\bR^{d+1}_{+})}
\leq N(1-e^{-\lambda})^{-(p-1)/p}\|f\|_{L_{p}(\bR^{d+1}_{+})},
$$
where the second inequality follows
from H\"older's inequality. 
\end{remark}

Similar improvement of integrability occurs for
the Green's function of $x_{t}$ rather than $(t,x_{t})$.
Notice  that 
$$
g_{\cA}(x):=\int_{0}^{\infty}G_{\cA}(t,x)\,dt
$$
satisfies
$$
EI_{\cA}\int_{0}^{\infty}e^{- \lambda s}
f( x_{\tau+s}-x_{\tau}) \,ds =\int_{\bR^{d }}f( x)g_{\cA}( x)\,dx 
$$
for any Borel nonnegative $f$ on $\bR^{d}$. For this reason
we call $g_{\cA}$ the Green's function of $x_{\tau+s}-x_{\tau}$ on
$\cA$.  

By using Remark \ref{remark 8.24.1} and Corollary \ref{corollary 10.4.1}
with $p=d$, $q=\infty$ we come to the following.
\begin{theorem}
                                       \label{theorem 10.4.1}
 We have
\begin{equation}
                                         \label{10.4.5}
 \|\Psi^{-1}_{\lambda} g_{\lambda}\|_{L_{d/(d-1) }(\bR^{d}) } 
\le N(d,\delta)\lambda^{-1/2}P(\cA ),
\end{equation}
where   $\Psi _{\lambda}( x)=\exp(- \sqrt{\lambda}
  |x| \bar\xi/16)$.
 
\end{theorem}

According to this theorem 
this Green's function is summable to the power $d/(d-1)$.
Again it turns out that this power can be increased.
If $B$ is an open ball in $\bR^{d}$ by $2B$ we denote
the concentric open ball of twice the radius of $B$.  

\begin{theorem}
                                  \label{theorem 10.4.3}
There exist  $d_{0}\in(1,d)$ and a constant    $N $, depending only
on $d,\delta$, 
 such that for any   ball  $B$ of radius $R\leq \kappa (\bar\xi)/(2\sqrt\lambda)$ 
 and $p\geq d_{0}  $, we have
\begin{equation}
                                          \label{10.4.6}
\| g_{\cA}\|_{L_{p/(p-1)}(B)}\leq N R^{-d /p }
\| g_{\cA}\|_{L_{1}(2B )} , 
\end{equation}
which is equivalently rewritten as
$$
\Big(\dashint_{B}g^{p/(p-1)}_{\cA}\,dx \Big)^{(p-1)/p}
\leq N\dashint_{2B }g_{\cA}\,dx.
$$
 
\end{theorem}

Proof. We again follow the idea in \cite{FS_84}.
Take a ball  $B$ of radius $R\leq \kappa (\bar\xi)/(2\sqrt\lambda)$ and
on the set $\{\tau<\infty\}$
define 
recursively  $\tau^{0}$ as the first  time
  when the process $ x_{\tau+t}-x_{\tau} $
hits $\bar B$, $\gamma^{0}$ as the first time after $\tau^{0}$
when this process leaves $2B$, $\tau^{n }$ as the first  time
after $\gamma^{n-1}$ when the process $ x_{\tau+t}-x_{\tau} $
hits $\bar B$, $\gamma^{n }$ as the first time after $\tau^{n }$
when this process leaves $2B$.
 
Then for any nonnegative Borel $f$ vanishing outside $B$
with $\|f\|_{L_{d }(B)}=1$ 
we have
$$ 
I:=
\int_{B}f ( x)  g_{\cA}( x)\,dx
=\sum_{n=0}^{\infty}EI_{\cA}e^{- \lambda \tau^{n} } 
E_{\cF_{\tau+\tau^{n}}} \int_{0}^{\gamma^{n}-\tau^{n}}e^{-\lambda t }f( x_{\tau+\tau^{n}+t}-x_{\tau})\,dt.
$$

Next we use   
\eqref{9.29.5}
 to see that the conditional expectations
above are less than $NR $. After that
we use   Corollary   
\ref{corollary 7.29.1} 
to get that   
$$
R^{2}I_{\tau^{n}<\infty}\leq N E_{\cF_{\tau^{n}}} \int_{\tau^{n}}^{\gamma^{n}}
e^{-\lambda(t-\tau^{n})} \,dt.
$$
Then we obtain  
$$
I
\leq N  R ^{-1}
\sum_{n=1}^{\infty}EI_{\cA} e^{-\lambda \tau^{n} }
\int_{\tau^{n}}^{\gamma^{n}}
e^{-\lambda(t-\tau^{n})} \,dt
=N  R ^{-1}
\sum_{n=1}^{\infty}EI_{\cA}  
 \int_{\tau^{n} }^{\gamma^{n} }e^{-\lambda t } \,dt
$$
$$
\leq N  R ^{-1}
EI_{\cA} \int_{0}^{\infty}e^{-\lambda t}I_{2B}( x_{t}) \,dt
=N  R ^{-1}\int_{2B}g_{\cA}( x)\,dx .
$$

The arbitrariness of $f$ 
implies that
$$
\Big(\dashint_{B}  g_{\cA} ^{d/(d-1)}(x)\,dx 
\Big)^{(d-1)/d}\leq N \dashint_{2B}  g_{\cA}(x)\,dx,
$$
and again it only remains to use Gehring's lemma
in case $p=d$.
For larger $p$ it suffices to use H\"older's inequality.
The theorem is proved.  \qed

By mimicking the proof of Theorem \ref{theorem 2.3.1}
one gets its ``elliptic'' counterpart.

\begin{theorem}
                                      \label{theorem 2.3.2}
For any $p\geq d_{0}$  we have
$$
\|g_{\cA}\|_{L_{p/(p-1)}(\bR^{d})}
 \leq N(\delta,d )
\lambda ^{d/(2p)-1}.
$$
In particular, for any Borel $f\geq0$ given on $\bR^{d}$ (and $\lambda\geq \kappa^{2}(\bar\xi)\rho_{b}^{-2}$),
on $\{\tau<\infty\}$ we have
\begin{equation}
                               \label{8.26.4}
E_{\cF_{\tau}}\int_{0}^{\infty}e^{-\lambda t}f(x_{\tau+t})\,dt
\leq N(\delta,d )
\lambda ^{d/(2d_{0})-1}\|f\|_{L_{d_{0}}(\bR^{d})}.
\end{equation}
\end{theorem}

\begin{remark}
                                           \label{remark 2.7.1}
Below by $d_{0}$ we denote {\em any\/}
constant in $(1,d)$ for which
both estimates \eqref{8.26.5}
from Theorem \ref{theorem 2.3.1} and
\eqref{8.26.4} from Theorem
 \ref{theorem 2.3.2} hold for any Borel $f\geq0$ and stopping time $\tau$. This $d_{0}$ may depend
on the  process $x_{t}$ and it is always
smaller than or equal to $ d_{0}(d,\delta)$,
which is the maximum of the $d_{0}$'s 
from Theorems \ref{theorem 2.3.1} and 
 \ref{theorem 2.3.2}

Observe that, as the simple example of $a^{ij}=\delta^{ij}$
and $b\equiv0$
shows, $d_{0}(d,\delta)>d/2$ and $d_{0}(d,1)$
can be taken to be as close to $d/2$  as we wish. We call $ d_{0}(d,\delta)$ {\em the Fabes-Stroock constant\/} because these authors discovered and proved
in \cite{FS_84}
its existence in terms of PDEs.  
The reason for the above ``any''  is that sometimes there are other ways to obtain estimates \eqref{8.26.5} and
\eqref{8.26.4} without using Gehring's lemma. 
\end{remark}

By virtually repeating the proof of
Theorem 3.4 of \cite{Kr_20_2} based on an interpolating argument,
  we come to the following.

\begin{theorem}
                         \label{theorem 8.27.1}
Suppose     the following condition $(d_{0},p,q)$:
\begin{equation}
                         \label{10.7.1}
 p,q\in(1,\infty] ,\quad
\nu:=1-\frac{d_{0}}{p }-\frac{1}{q }\geq 0
\end{equation}
holds. 
Then for any Borel $f\geq0$ given on $\bR^{d+1}$ 
(and recall that $\lambda\geq \kappa^{2}(\bar\xi)\rho_{b}^{-2}$)
we have on $\{\tau<\infty\}$ that
\begin{equation}
                                 \label{8.27.01}
E_{\cF_{\tau}}\int_{0}^{\infty}e^{-\lambda t}f(t,x_{\tau+t})\,dt
\leq N(d,\delta )
\lambda^{ (1/2)(d/p+2/q)-1}\|f\|_{L_{p,q}}.
\end{equation}

\end{theorem}

By using the same argument as in Remark \ref{remark 10.15.1}
we can replace $(t,x_{\tau+t})$ with $(\tau+t,x_{\tau+t})$ in
\eqref{8.27.01} and then we see that Assumption
\ref{assumption 8.24.1} is satisfied for $p,q$ as in
 $(d_{0},p,q)$  \eqref{10.7.1}, $\ell=N(d,\delta )$ and $\alpha=
1-(1/2)(d/p+2/q)$. Then Theorem \ref{theorem 8.22.1}
is valid, which yields the following.

\begin{theorem}
                                  \label{theorem 8.30.1}
 
Assume  $(d_{0},p,q)$  \eqref{10.7.1}.
Then for any Borel $f\geq0$ given on $\bR^{d+1}$ 
and $\lambda\geq \kappa^{2}(\bar\xi)\rho_{b}^{-2}$ 
we have on the set $\{\tau<\infty\}$ that
\begin{equation}
                                 \label{8.27.1}
E_{\cF_{\tau}}\int_{0}^{\infty}e^{- \lambda s}
f(s,x_{\tau+s}-x_{\tau}) \,ds
\leq N(\delta,d )
\lambda^{ -\chi}\|\Psi_{\lambda}f\|_{L_{p,q}(\bR^{d+1}_{+})}.
\end{equation}
where $\Psi _{\lambda}(t,x)=\exp(- 
\sqrt{\lambda} (|x|+ \sqrt t)\bar\xi/16)$, $\chi=1-(1/2)(d/p+2/q)$.
In particular, if $f$ is independent of $t$, $p\geq d_{0}$,
   $q=\infty$, and $\bar \Psi _{\lambda}( x)=\exp(- 
\sqrt{\lambda}  |x| \bar\xi/16)$,
$$
E_{\cF_{\tau}}\int_{0}^{\infty}e^{- \lambda t}
f( x_{\tau+t}-x_{\tau}) \,dt\leq 
N\lambda ^{-1+d/(2p)}\|\bar \Psi_{\lambda}^{d_{0}/p} f\|_
{L_{p   } },
$$
 
\end{theorem}

\begin{theorem}
                           \label{theorem 9.7.1}
Assume that  $(d_{0},p,q)$  \eqref{10.7.1} holds.
Then

(i)
  for any
$n=1,2,...$, nonnegative Borel $f$ on $\bR^{d+1}_{+}$, 
  and
 $T\leq \kappa^{-2}(\bar\xi)\rho_{b}^{2}$  we have
on $\{\tau<\infty\}$ that
\begin{equation}
                                          \label{9.7.1}
E_{\cF_{\tau}}\Big[\int_{0}^{T}  
f(t,x_{\tau+t}-x_{\tau})\,dt\Big]^{n}\leq n!N^{n} (d,\delta)
T^{n\chi }\| \Psi^{(1-\nu)/n} _{1/T}
f\|^{n}_{L_{p,q}(\bR^{d+1}_{+}) },
\end{equation}

(ii)  for any
  nonnegative Borel $f$ on $\bR^{d+1}_{+}$,
$R\leq 1$,  and
 $T\geq \kappa^{-2}(\bar\xi)\rho_{b}^{2}$  we have
on $\{\tau<\infty\}$ that
\begin{equation}
                               \label{9.7.10}
I:=E_{\cF_{\tau}} \int_{0}^{T}  
f(t,x_{\tau+t} )\,dt \leq N(d,\delta,\rho_{b}) TR^{-2-d}\sup_{C\in\cC_{R}}
 \|  
f\| _{L_{p,q}(C) } .  
\end{equation}

\end{theorem}

Proof. To prove   (i) we proceeds by induction
on $n$. The induction hypothesis is that for   $\kappa\in[0,1/n]$ and any $\bR^{d+1}_{+}$-valued $\cF_{\tau}$-measurable $(\gamma,\xi)$
$$
E_{\cF_{\tau}}\Big[\int_{0}^{T}  
f(\gamma+t, x_{\tau+t}-x_{\tau}+\xi )\,dt\Big]^{n}
$$
\begin{equation}
                                   \label{10.28.1}
\leq n!N^{n} 
T^{n\chi } \Psi_{1/T}^{( \nu-1)\kappa n}
(\gamma,\xi)\| \Psi^{( 1-\nu)\kappa} _{1/T}
f\|^{n}_{L_{p,q}(\bR^{d+1}_{+}) }.
\end{equation}

Denote $\Delta_{t}x_{\tau}=x_{\tau+t}-x_{\tau}$ and observe that for $s\geq \tau+t$ we have $\Delta_{s}x_{\tau}
=\Delta_{s-t}x_{\tau+t}+\Delta_{t}x_{\tau}$.
If the hypothesis  holds true for some $n\geq1$,
then  by observing that
$$
J :=E_{\cF_{\tau}}\Big[\int_{0}^{T}  
f(\gamma+t,\Delta_{t}x_{\tau}+\xi)\,dt\Big]^{n+1}
$$
$$
=(n+1)E_{\cF_{\tau}}\int_{0}^{T}f(\gamma+t,\Delta_{t}x_{\tau}+\xi)
\Big(E_{ \cF_{\tau+t} } \Big[\int_{ t}^{T} 
f(\gamma+s,\Delta_{s}x_{\tau}+\xi)\,ds\Big]^{n} \,ds\Big)dt,
$$
we see that, for any $\kappa\in[0,1/n]$,
\begin{equation}
                                        \label{8.31.1}
J \leq (n+1)!N^{n}T^{n\chi]}\|\Psi^{(1-\nu)\kappa}_{1/T}
f\|^{n}_{L_{p}(\bR^{d})}
 E_{\cF_{\tau}}\int_{0}^{T}\Psi^{(\nu-1)\kappa n}_{1/T}f(\gamma+t,\Delta_{t}x_{\tau}+\xi)\,dt.
\end{equation}

  We have for any $\lambda>0$
$$
 E_{\cF_{\tau}}\int_{0}^{T}\Psi^{(\nu-1)\kappa n}_{1/T}f(\gamma+t,\Delta_{t}x_{\tau}+\xi)\,dt
$$
$$
\leq e^{\lambda T}
E_{\cF_{\tau}}\int_{0}^{\infty}e^{-\lambda t}\Psi^{(\nu-1)\kappa n}_{1/T}f(\gamma+t,\Delta_{t}x_{\tau}+\xi)\,dt.
$$
 Here the last term,
owing to Theorem \ref{theorem 8.30.1}, for $\lambda=1/T$ 
and $\mu\in[0,1]$ is
dominated by
$$
N(\delta,d ) T^{\chi} \|\Psi^{(\nu-1)\kappa n}_{1/T}f(\gamma+\cdot,\xi+\cdot)\Psi^{\mu}_{1/T} \|_{L_{p,q}(\bR^{d+1}_{+})}
$$
$$
\leq N(\delta,d ) T^{\chi}\Psi^{-\mu}(\gamma,\xi) \|\Psi^{(\nu-1)\kappa n+\mu}_{1/T}f(\gamma+\cdot,\xi+\cdot)  \|_{L_{p,q}(\bR^{d+1}_{+})},
$$
where the last inequality is due to the fact that
$\Psi _{\lambda}(s,y)\leq \Psi  _{\lambda}(t+s,x+y)
\Psi^{-1} _{\lambda}(t, x)$.
For   $\kappa=\mu$
we get \eqref{9.7.1} with $n=1$ by replacing 
$\Psi^{(\nu-1)\kappa n }_{1/T}f$ by $f$,
which justifies the start of the induction.
For $\mu=(1-\nu)\kappa(n+1)$, $\kappa\in[0,1/(n+1)]$, we have
$\Psi^{(\nu-1)\kappa n+\mu}_{1/T}=
\Psi^{(1-\nu)\kappa   }_{1/T}$ and this along with 
\eqref{8.31.1} show  that our hypothesis holds true
also for $n+1$.
This proves \eqref{9.7.1}.

While proving \eqref{9.7.10} we may assume that  
$R=1$ (see Remark \ref{remark 2.29.1} below) and that
  $T=k \beta$, where $k\geq 1$ is an integer
and $\beta=\kappa^{-2}(\bar\xi)\rho_{b}^{2}$. Then 
first consider the case of $\nu=0$. Note that
owing to \eqref{9.7.1}  
\begin{equation}
                                    \label{12.18.03}
E_{\cF_{\tau}}\int_{0}^{\beta}f(t,x_{\tau+t})
\,dt\leq N(d,\delta)\rho_{b}^{(2d_{0}-d)/p}\|f(\cdot,\cdot+x_{\tau})\Psi_{1/\beta}
 \|_{L_{p,q}(\bR^{d+1}_{+})} .
\end{equation}

Let $\cZ=\{0,1,2,...\}\times \bZ^{d}$ and for
$z=(z_{1},z_{2})\in\cZ$ let $C^{z}=C_{1}(z)$. Observe that
on $C^{z}$ we have  
$$
\Psi_{1/\beta} \leq  \exp(- 
2\mu(|z_{2}|+ \sqrt z_{1} ) ),
$$
where $2\mu=\beta^{-1/2} \bar\xi/16$.
Furthermore, for each $z\in\cZ$
$$
\|f(\cdot,\cdot+x_{\tau})\|_{L_{p,q}(C^{z})}\leq \sup_{C\in\cC_{1}}
 \|  
f\| _{L_{p,q}(C) }.
$$
Therefore, by noting that
$ f(\cdot,\cdot+x_{\tau})\Psi_{1/\beta}  \leq \sum_{\cZ}  f(\cdot,\cdot+x_{\tau} )\Psi_{1/\beta} 
I_{C^{z}}$ and using Minkowski's inequality we get
that the norm in \eqref{12.18.03} is dominated by
$$
\sup_{C\in\cC_{1}}
 \|  
f\| _{L_{p,q}(C) }\sum_{\cZ}
 \exp(- 2
\mu(|z_{2}|+ \sqrt z_{1} ) ).
$$
By majorating the last sum by an integral
  we obtain
that it is dominated by
$$
\int_{0}^{\infty}\int_{\bR^{d}}e^{-2\mu (|x|+ \sqrt t-3)_{+} }\,dxdt\leq N+
\int_{0}^{\infty}\int_{\bR^{d}}e^{-2\mu (|x|+ \sqrt t-3)_{+} }I_{|x|+\sqrt t>6}\,dxdt
$$
$$
\leq N+
\int_{0}^{\infty}\int_{\bR^{d}}e^{- \mu (|x|+ \sqrt t )  }I_{|x|+\sqrt t>6}\,dxdt
$$
$$
\leq N+
\int_{0}^{\infty}\int_{\bR^{d}}e^{- \mu (|x|+ \sqrt t )  } \,dxdt=N+N\mu^{-d-2}.
$$
Hence, for $n=0$
$$
E_{\cF_{\tau+n\beta}}\int_{n\beta}^{(n+1)\beta}f(t,x_{\tau+t})
\,dt\leq N(d,\delta)\rho_{b}^{(2d_{0}-d)/p}
(1+\rho_{b}^{d+2})
\sup_{C\in\cC_{1}}\|
f\| _{L_{p,q}(C) }.
$$
Clearly, this also holds for any $n=1,2,...$ and since
$T=k\beta=k\kappa^{-2}(\bar\xi)\rho_{b}^{2}$,
$$
I \leq N(d,\delta)T\rho_{b}^{-2}\rho_{b}^{(2d_{0}-d)/p}
(1+\rho_{b}^{d+2})
\sup_{C\in\cC_{1}}\|
f\| _{L_{p,q}(C) } 
$$
and this  
 proves \eqref{9.7.10} if $\nu=0$.

If $\nu=1$ ($p=q=\infty$), \eqref{9.7.10} is obvious,
and if $\nu<1$, by the above, \eqref{9.7.10}
holds with $(1-\nu)(p,q)$ in place of $(p,q)$,
which yields \eqref{9.7.10} as is after using 
H\"older's inequality.
The theorem is 
proved. \qed

\begin{remark}
                     \label{remark 2.29.1}
If $R\leq 1$, the term $\sup_{C\in\cC_{1}}
 \|  
f\| _{L_{p,q}(C) }$ in \eqref{9.7.10}
can be replaced
with 
$
R^{-2-d}\sup_{C\in\cC_{R}}
 \|  
f\| _{L_{p,q}(C) }.
$ 

Indeed, by simple inspection one proves that
for any $R\geq 1$, $C\in\cC_{R}$ and $C'$,
defined as the union of $2C$ and its reflection
in its lower base,
$$
 I_{ C}( t,  x)\leq N\int_{2 C}I_{C_{1}}( t- s,  x-  y)\,
dsdy,
$$
where $N=N(d)$. Dilations show that, for any
$R\leq 1$, $r\geq R$, and $C\in C_{r}$ we have
$$
 I_{ C}( t,  x)\leq NR^{-2-d}\int_{2 C}I_{C_{R}}( t- s,  x-  y)\,
dsdy.
$$
 It follows
for $C\in \cC_{1}$ and $R\leq1$ that
$$
 |f|I_{ C}\leq NR^{-2-d}\int_{2 C}|f|
I_{C_{R}(s,y)}\,
dsdy,
$$
$$
\|f\|_{L_{p,q}(C)}\leq NR^{-2-d}\int_{2C}
\sup_{C\in\cC_{R}}\|
f\| _{L_{p,q}(C) }\,dsdy=
NR^{-2-d}
\sup_{C\in\cC_{R}}\|
f\| _{L_{p,q}(C) }.
$$
\end{remark}

Next theorem improves 
estimate \eqref{9.29.2} in what concerns
the restrictions on $p,q$.
\begin{theorem}
                          \label{theorem 9.5.1}
Assume that  $(d_{0},p,q)$  \eqref{10.7.1} holds.
 Then 
for any $R\leq \kappa^{-1} (\bar\xi)\rho_{b} $, $\cF_{\tau}$-measurable $\bR^{d}$-valued $y$,  and Borel nonnegative $f$ 
given on $\bR^{d+1}$,
we have on $\{\tau<\infty\}$ that
\begin{equation}
                                \label{9.5.4}
E_{\cF_{\tau}}\int_{0}^{\theta_{\tau}\tau_{ R}(y) }f( t, x_{\tau+t}-x_{\tau}-y)\,dt\leq
N(d,\delta)R^{2}\dashnorm f\|_{L_{p,q}(C_{R})},
\end{equation}
\begin{equation}
                             \label{9.3.1}
E_{\cF_{\tau}}\int_{0}^{\theta_{\tau}\tau'_{ R}(y)  }f( t, x_{\tau+t}-x_{\tau})\,dt\leq
N(d,\delta)R^{2}
\sup_{C\in\cC_{R}}
 \dashnorm  
f\| _{L_{p,q}(C) }.
\end{equation}

\end{theorem}

Proof.  Since $\theta_{\tau}\tau_{ R}(y) \leq R^{2}$, the left-hand side of
\eqref{9.5.4} is smaller than
$$
e^{\lambda R^{2}}E_{\cF_{\tau}}\int_{0}^{\infty}e^{-\lambda t}
I_{C_{R}}f(t,x_{\tau+t}-x_{\tau}-y)\,dt
$$
for any $\lambda>0$. For $\lambda=R^{-2}$ we have
$\lambda\geq \kappa^{2}(\bar\xi)\rho_{b}^{-2}$
and \eqref{9.5.4} follows from \eqref{8.27.1}.

To prove \eqref{9.3.1}, it suffices to note that \eqref{9.5.4} remains valid if its right-hand side is replaced with that of \eqref{9.3.1}, and then repeat the same argument
as in the proof of \eqref{1.3.3}.
 The theorem
is proved. \qed

Here is a key to finding analytic
conditions insuring  that $\bar b_{\rho_{b}}\leq m_{b}$.

\begin{corollary}
                                       \label{corollary 3.26.1}
Assume that there exists 
  functions $b_{i}(t,x)\geq0$,
$i=1,...,k$, on $\bR^{d+1}$ such that
$|b_{t}|\leq (b_{1}+...+b_{k})(t,x_{t})$ for all $(\omega,t)$. Take some
$p_{i},q_{i}$ satisfying 
$(d_{0},p,q)$  \eqref{10.7.1}. Suppose that
there is
a constant $\hat b\in(0,\infty)$
such that, for any $\rho\leq \hat\rho_{b}
:=\kappa^{-1} (\bar\xi)\rho_{b}$ and $C\in\cC_{\rho}$
\begin{equation}
                                                \label{3.26.2}
\sum_{i=1}^{k}\dashnorm b_{i}\|_{L_{p_{i},q_{i}}(C)}\leq\hat b \rho^{ -1}.
\end{equation}
Then $\bar b_{ \hat\rho_{b}}\leq  N
(d,\delta)\hat b$. 

\end{corollary}

Indeed, \eqref{9.5.4} implies that
$$
E_{\cF_{t}}\int_{0}^{\theta_{t}\tau_{ R}(x) }(b_{1}+...+b_{k})( t+s, x_{t+s}  )\,ds
$$
$$
=N(d,\delta)R^{2}\sum_{i=1}^{k}\dashnorm b_{i}I_{C_{R}(t,x_{t}+x) }\|_{L_{p_{i},q_{i}} }
\leq N(d,\delta)\hat b R .
$$

 \begin{remark}
                                 \label{remark 10.18.1}
In light of Corollary \ref{corollary 3.26.1} it is tempting
to claim that if \eqref{3.26.2} holds and $N(d,\delta)\check b
\leq m_{b}$, then our main Assumption   $(\bar  b_{\rho_{b}})$  \ref{assumption 8.19.2}
is satisfied.  However this is a vicious circle:
Corollary \ref{corollary 3.26.1} was obtain on the basis
of Assumption   $(\bar  b_{\rho_{b}})$  \ref{assumption 8.19.2}. Nevertheless,
it turns out in the case of stochastic equations that
if condition \eqref{3.26.2} is satisfied and 
$N(d,\delta)\check b
\leq m_{b}$, then there exists at least one solution
of the equation, for which Assumption   $(\bar  b_{\rho_{b}})$  \ref{assumption 8.19.2}
is satisfied.

\end{remark}

One also has an estimate similar to 
\eqref{9.3.1} for $R>\kappa^{-1}(\bar \xi)\rho_{b}$, however, with not so sharp control of the constants.

\begin{theorem} 
                      \label{theorem 3.27.20}

Assume that  $(d_{0},p,q)$  \eqref{10.7.1} holds.   
Then
for   any $R\in(0,\infty)$,   $\cF_{\tau}$-measurable $\bR^{d}$-valued $y$,
and Borel $f\geq0$ on $\{\tau<\infty\}$ we have
\begin{equation}
                                              \label{3.21.80}
E_{\cF_{\tau}}\int_{0}^{\theta_{\tau}\tau'_{ R}(y)  }f( t, x_{\tau+t}-x_{\tau})\,dt\leq
\hat N   \sup_{C\in\cC_{1}}
 \|  
f\| _{L_{p,q}(C) } , 
\end{equation}
where $\hat N$ depends only on 
$d,\delta ,\rho_{b}$, and $R $.
\end{theorem}

Proof.   By Corollary \ref{corollary 10.17.1} we have
$P_{\cF_{\tau}}(\theta_{\tau}\tau'_{ R}(y) >T)\leq Ne^{-T/N}$
for all $T$ with $N=N(d,\delta,\rho_{b},R)$ and there exists
$T=T(d,\delta,\rho_{b},R)$ such that the right-hand side is less than $1/2$.
This shows, by the same argument
as in the proof of \eqref{1.3.3},     
that to prove the current theorem 
it suffices to prove that  
\begin{equation}
                                              \label{3.1.4}
E_{\cF_{\tau}}\int_{0}^{T\wedge\theta_{\tau}\tau'_{ R}(y)   }f( t, x_{\tau+t}-x_{\tau} )\,dt\leq
\hat N  \sup_{C\in\cC_{1}}
 \|  
f\| _{L_{p,q}(C) } ,  
\end{equation}
where $\hat N$ depends only on 
$d,\delta ,\rho_{b}$, and $R$. Here the left-hand side
is less than
$$
E_{\cF_{\tau}}\int_{0}^{T    }f( t, x_{\tau+t}-x_{\tau} )\,dt,
$$
so that \eqref{3.1.4} follows from \eqref{9.7.10}.
The theorem is proved. \qed

Theorem \ref{theorem 9.7.1} allows us to prove It\^o's formula
for functions $u\in W^{1,2}_{p,q}(Q)$, where $Q$
is a domain in $\bR^{d+1}$ and 
$$
W^{1,2}_{p,q }(Q)=\{v: v, \partial_{t}v,
 Dv,  D^{2}v\in L_{p,q }(Q) \}
$$
with norm introduced in a natural way.
Before, the formula was known only for
(smooth (It\^o) and) $W^{1,2}_{d+1}$-functions
and processes with bounded drifts or
for $W^{2}_{d_{0}}$-functions in case the drift
of the process is dominated by $h(x_{t})$
with $h\in L_{d}$ (see \cite{Kr_19_1}).

The following extends Theorem 2.10.1  of \cite{Kr_77}
to functions with lower summability of the derivatives
and to the spaces with mixed norms.

\begin{theorem}
                                 \label{theorem 10.15.1}
Assume that  $(d_{0},p,q)$  \eqref{10.7.1} holds with    $p<\infty$,
$q<\infty$. Let
$Q$ be a bounded domain in $\bR^{d+1}$, $0\in Q$,
$b$ be {\em bounded\/}, and  $u\in W^{1,2}_{p,q}(Q)\cap C(\bar Q)$. Then,
for $\tau$ defined as the first exit time of $(t,x_{t})$
from $Q$ with  probability one for all $t\geq0$,
$$
u(t\wedge\tau,x_{t\wedge\tau})
=u(0,0)+\int_{0}^{t\wedge\tau}D_{i}u(s,x_{s})\sigma^{ik}_{s}\,dw^{k}_{s}
$$
\begin{equation}
                                      \label{10.15.01}
+\int_{0}^{t\wedge\tau}[
\partial_{t}u(s,x_{s})+ a^{ij}_{s}D_{ij}u(s,x_{s})
+b^{i}_{s}D_{i}u(s,x_{s})]\,ds
\end{equation}
and the stochastic integral above is a square-integrable
martingale.
\end{theorem}

Proof. First assume that $u$ is smooth and its derivatives
 are bounded. Then  
\eqref{10.15.01} holds by It\^o's formula and, moreover,
by denoting $\tau^{n}= n\wedge\tau$
for any $n\geq0$ we have
$$
E \int_{ \tau^{n}}^{\tau^{n+1}}|Du(s,x_{s})|^{2}\,ds
\leq NE\Big(\int_{\tau^{n}}^{ \tau^{n+1}
}D_{i}u(s,x_{s})\sigma^{ik}_{s}\,dw^{k}_{s}\Big)^{2}
$$
$$
=NE\Big(u( \tau^{n+1},x_{ \tau^{n+1}})-
u( \tau^{n },x_{ \tau^{n }})
$$
$$-
\int_{ \tau^{n}}^{\tau^{n+1}}[
\partial_{t}u(s,x_{s})+ a^{ij}_{s}D_{ij}u(s,x_{s})
+b^{i}_{s}D_{i}u(s,x_{s})]\,ds\Big)^{2}
$$
$$
\leq N\sup_{\bar Q}|u|
+NE\Big(\int_{ \tau^{n}}^{\tau^{n+1}}I_{Q}
\big(|\partial_{t}u|+|Du|+|D^{2}u|\big)(s,x_{s})\,ds\Big)^{2}.
$$
Since $Q$ is bounded, $\tau$ is bounded as well
and
in light of Theorem \ref{theorem 9.7.1} we conclude that
\begin{equation}
                                      \label{10.15.2}
E \int_{0}^{\tau}|Du(s,x_{s})|^{2}\,ds\leq N
\sup_{\bar Q}|u|+N \| \partial_{t}u,
 Du, D^{2}u\|_{L_{p,q}(Q)},
\end{equation}
where $N$ are independent of $u$ and $Q$
as long as the size of $Q$ in the $t$-direction
is under control.
Owing to Fatou's theorem,
this estimate is also true for those $u\in W^{1,2}_{p,q}(Q)
\cap C(\bar Q)$
that can be approximated uniformly and in the 
$W^{1,2}_{p,q}(Q)$-norm by smooth functions with bounded
derivatives (recall that $p<\infty$,
$q<\infty$). For our $u$
there is no guarantee that such approximation is possible.
However, mollifiers do such approximations
 in any subdomain $Q'\subset \bar Q'\subset
Q$. Hence, \eqref{10.15.2} holds for our $u$ if we replace $Q$
by $Q'$
(containing $(0,0)$). Setting $Q' \uparrow Q$  proves \eqref{10.15.2}
in the generals case and proves the last assertion
of the theorem.

After that \eqref{10.15.01} with $Q'$ in place of $Q$
is proved by routine approximation of $u$ by smooth 
functions. Setting $Q' \uparrow Q$ finally proves \eqref{10.15.01}.
The theorem is proved.\qed

Here is an ``elliptic'' version 
based on Theorem  \ref{theorem 9.7.1} with $q=\infty$.
 For $p\geq d$ 
Theorem \ref{theorem 12.23.1} can be found in \cite{Kr_77}.
\begin{theorem}
                          \label{theorem 12.23.1}
Assume that  $p\in[d_{0},\infty)$.   Let
$G$ be a bounded domain in $\bR^{d}$, $0\in Q$,
$b$ be {\em bounded\/}, and  $u\in W^{ 2}_{p }(G)\cap C(\bar G)$ ($u$ is independent of $t$). Then,
for $\tau$ defined as the first exit time of $ x_{t} $
from $G$ with  probability one for all $t\geq0$,
$$
u( x_{t\wedge\tau})
=u(0 )+\int_{0}^{t\wedge\tau}D_{i}u( x_{s})\sigma^{ik}_{s} \,dw^{k}_{s}
+\int_{0}^{t\wedge\tau}[
 a^{ij}_{s}D_{ij}u( x_{s})
+b^{i}_{s}D_{i}u( x_{s})]\,ds 
$$
and the stochastic integral above is a square-integrable
martingale.
\end{theorem}


\begin{thebibliography}{mm}

\bibitem{Ca_95} X. Cabr\'e, 
{\em On the Alexandroff-Bakelman-Pucci estimate and the reversed
H\"older inequality for solutions of elliptic 
and parabolic equations\/}, Comm. Pure
Appl. Math., 48 (1995), 539--570.

\bibitem{CKS_00} M. G. Crandall, M. Kocan, and A. \'Swi{\c e}ch, {\em
$L^p$-theory for fully nonlinear uniformly parabolic equations\/},
 Comm. Partial Differential Equations, Vol. 25  (2000),
 No. 11-12, 1997--2053.


\bibitem{DK_21}  Hongjie Dong and N.V. Krylov,
{\em
Aleksandrov's estimates for elliptic equations with drift in a Morrey
  spaces containing $L_{d}$\/},  http://arxiv.org/abs/2103.03955

\bibitem{FS_84} E.B. Fabes and D.W. Stroock,
{\em The $L^{p}$-integrability of Green's functions
and fundamental solutions for elliptic
and parabolic equations\/}, Duke Math. J., Vol. 51
(1984), No. 4, 997--1016.

\bibitem{Fo_98} K. Fok, {\em A nonlinear Fabes-Stroock result\/},
Comm. Partial Differential Equations, 23 (1998), No. 5-6, 967--983.

\bibitem{GS_82} M.Giaquinta and M. Struwe,
{\em On the partial regularity of weak solutions 
of nonlinear parabolic systems\/}, 
Mathematische Zeitschrift, Vol. 179 (1982), 437-451. 

\bibitem{Kr_74} N.V. Krylov, 
{\em  Some estimates for the density of  distribution
  of a
 stochastic
 integral\/}, Izvestiya Akademii Nauk SSSR, seriya matematicheskaya,
Vol. 38 (1974), No. 1,  228--248 in Russian; English translation
in Math. USSR
Izvestija,   Vol. 8  (1974),  No. 1,   233--254.

\bibitem{Kr_77}  N.V. Krylov,  ``Controlled diffusion processes'',
Nauka, Moscow,  1977 in Russian; English  transl. 
   Springer,
1980.

\bibitem{Kr_86} N.V. Krylov,
  {\em  On  estimates of the maximum of a solution of
a parabolic equation and  estimates of the distribution of a
semimartingale}, Matematicheski Sbornik, Vol. 130, No. 2  (1986),  207--221
in Russian, English translation is
 Math. USSR Sbornik,  Vol. 58 (1987), No. 1,  207--222.


\bibitem{Kr_19_1} N.V. Krylov, {\em
On stochastic equations with drift in $L_{d}$\/},
 Stoch. Proc. Appl., Vol. 138 (2021), 1--25.

\bibitem{Kr_20_2} N.V. Krylov, {\em
On time inhomogeneous stochastic It\^o equations 
with drift in $L_{d+1}$\/},  
Ukrains'kyi Matematychnyi Zhurnal, Vol. 72 (2020), No. 9,  1232--1253, 
reprinted in
 Ukrainian Math. J. 72 (2021), no. 9, 1420--1444. 

\bibitem{4} N.V. Krylov, {\em On   diffusion
 processes with drift in a Morrey class containing  $L_{d+2}$\/}, J. Dyn. Diff. Equat. (2021),
 https://doi.org/10.1007/s10884-021-10099-x.

\bibitem{2022_2} N.V. Krylov, {\em
Some properties of
solutions of  It\^o equations 
with drift in $L_{d+1} $\/},
  Stoch. Proc. Appl., Vol. 147 (2022), 363--387.

\bibitem{Na_15} A.I. Nazarov, 
{\em
Interpolation of linear spaces and estimates for 
the maximum of a solution for parabolic equations\/}, Partial differential
equations, Akad.  Nauk SSSR Sibirsk. Otdel., Inst. Mat., No\-vo\-si\-birsk,
1987, 50--72 in Russian; translated into English 
as  {\em
On the maximum principle for parabolic equations
with unbounded coefficients\/},
https:// arxiv.org/abs/1507.05232

\bibitem{NU_85} A.I. Nazarov and N.N. Ural'tseva, {\em
Convex-monotone hulls and an estimate of the maximum of the solution of
a parabolic equation\/}, Boundary value problems of
mathematical physics and related problems 
in the theory of functions,
No. 17,  Zap. Nauchn. Sem. Leningrad. Otdel. Mat. Inst. Steklov. 
(LOMI) Vol. 147 
(1985), 95--109, in Russian, English translation in
 Journal of Soviet
Mathematics Vol. 37 (1987), 851--859.
\end{thebibliography}
\end{document}